\documentclass[11pt]{article}
\usepackage[utf8]{inputenc}
\usepackage{textcomp}

\usepackage{graphicx}
\usepackage{amsmath,amsfonts,latexsym,cite}
\usepackage[version=4]{mhchem}
\usepackage{siunitx}
\usepackage{longtable,tabularx}
\usepackage{hyperref,amsfonts,latexsym,url,color,subfig,enumerate,graphics}
\usepackage{mathtools}
\usepackage{wasysym}

\title{\bf Minimum-Fuel Earth-Based Orbit Transfers Using Multiple-Domain Adaptive Radau Collocation}

\author{Brittanny V.~Holden\footnote{PhD. Student, Department of Mechanical and Aerospace Engineering, brittannyholden@gmail.com.} \\ Anil V.~Rao \footnote{Professor, Department of Mechanical and Aerospace Engineering, anilvrao@ufl.edu, Associate Fellow AIAA (Corresponding Author).} \vspace{12pt} \\ {\em University of Florida} \\ {\em Gainesville, FL 32611-6250}}
\date{}

\begin{document}

\maketitle

\begin{abstract}
  A numerical optimization study of minimum-fuel Earth-based orbital transfers from low-Earth orbit (LEO) to either medium-Earth orbit (MEO), high-Earth orbit (HEO), or geostationary orbit (GEO), is performed.  Various values of maximum allowable thrust acceleration are considered for each type of transfer (LEO-to-MEO, LEO-to-HEO, or LEO-to-GEO).  A key aspect of the study performed in this paper is that the optimal thrusting structure is not assumed to be known a priori, but is determined as part of the solution process.  In order to determine the optimal thrusting structure, a recently developed bang-bang and singular optimal control (BBSOC) method is employed together with multiple-domain Legendre-Gauss-Radau quadrature collocation.  Key results obtained in this study include not only the number of switches in the optimized thrust, but also the total impulse. Furthermore, it is found that, as the maximum allowable thrust acceleration decreases, the total impulse is less than the total impulse obtained from a previous study where a burn-coast-burn thrusting structure was assumed a priori.  For each type of transfer a particular value of maximum allowable thrust acceleration is chosen to highlight in more detail the key features of the optimal solutions.  This study provides improved results over previous studies and provides improved insight into the optimal thrusting structure required in order to accomplish each type of orbital transfer using the least amount of fuel.  
\end{abstract}

\renewcommand{\baselinestretch}{1}
\normalsize
\normalfont

\section*{Nomenclature}


{\renewcommand\arraystretch{1.0}
\noindent\begin{longtable*}{@{}l @{\quad=\quad} l@{}}
$a$ & semi-major axis \\ 
$A_{T}$ & number of thrust arcs \\ 
$\textrm{AU}$ & acceleration unit \\
$c_p$ & number of initial collocation points in each mesh interval \\ 
$\textrm{DU}$ & length unit \\
$e$ & eccentricity \\
$e_D$ & desired eccentricity of terminal orbit \\
$e_f$ & eccentricity of terminal orbit \\
$e_0$ & eccentricity of initial orbit \\
$f$ & second modified equinoctial orbital element \\
$\textrm{FU}$ & force unit \\
$g$ & third modified equinoctial orbital element \\
$g_0$ & standard gravitational acceleration \\
$h$ & fourth modified equinoctial orbital element \\ 
$i$ & orbital inclination \\ 
$i_D$ & desired inclination of terminal orbit \\
$i_f$ & inclination of terminal orbit \\
$i_0$ & inclination of initial orbit \\
$I_{sp}$ & specific impulse \\
$J$ & cost functional \\
$k$ & fifth modified equinoctial orbital element \\ 
$L$ & true longitude (sixth modified equinoctial orbital element) \\
$N$ & total revolutions around Earth \\
$M$ & number of initial mesh intervals \\
$m$ & mass of spacecraft \\
$m_0$ & initial mass of spacecraft \\
$\textrm{MU}$ & mass unit \\
$p$ & semi-parameter (first modified equinoctial orbital element) \\
$p_D$ & desired semi-parameter of terminal orbit \\
$p_f$ & semi-parameter of terminal orbit \\
$p_0$ & semi-parameter of initial orbit \\
$R_E$ & radius of Earth \\
$s_0$ & maximum allowable thrust acceleration \\
$t$ & time \\
$t_0$ & initial time \\ 
$t_f$ & terminal time \\ 
$t_{T}$ & total time thrusting \\
$T$ & thrust magnitude \\ 
$\textrm{TU}$ & time unit \\
$T_{\max}$ & maximum thrust magnitude \\
$\textrm{VU}$ & speed unit \\
$\textbf{u}$ & thrust direction \\
$u_n$ & normal component of thrust direction \\  
$u_r$ & radial component of thrust direction \\
$u_t$ & transverse component of thrust direction \\  
$\Delta_n$ & normal non-two-body perturbations of the spacecraft \\
$\Delta_r$ & radial non-two-body perturbations of the spacecraft \\
$\Delta_t$ & transverse non-two-body perturbations of the spacecraft \\
$\Delta V$ & the total impulse of the spacecraft \\
$\eta$ & threshold to detect the relative size of jumps in the control \\
$\mu_E$ & gravitational parameter of Earth \\
$\nu$ & true anomaly \\
$\omega$ & argument of periapsis \\ 
$\Omega$ & longitude of the ascending node \\ 
$\Omega_0$ & longitude of the ascending node of the initial orbit \\ 
\end{longtable*}}

\renewcommand{\baselinestretch}{2}
\normalsize
\normalfont

\section{Introduction}

A topic of great interest to the space community is the optimization of planet-based orbital transfers using advanced propulsion systems.  In the earlier days of space travel, high-thrust chemical propulsion was utilized.  More recently, however, low-thrust propulsion systems using either solar or electric propulsion have been used because such propulsion systems consume significantly less fuel compared with chemical propulsion, thereby extending mission life and reducing the launch costs because the overall mass of the payload is less relative to spacecraft that carry chemical propulsion resources.  This paper focuses on orbital transfers between Earth orbits using low-thrust propulsion.  

A significant amount of research has been conducted previously on the design of low-thrust orbital transfer missions~\cite{Morante2021, Rauwolf1996, Coverstone2000, Dachwald2005, Starek2012, Betts2014, Yang2009, JiangBaoyinLi2012, Gergaud2006, Daero2021, KlueverOleson1998, Betts2000, ScheelConway1994, SchubertRao2012, Graham2015, Herman2002}. Reference~\cite{Morante2021} primarily focused on solving hybrid-optimal control problems by conducting a survey on approaches for low-thrust trajectory optimization. Through substantial research of both traditional and new methods, it was determined that most optimizers that focus on low-thrust trajectories cannot include mission-planning easily, search over multi-objective design spaces, and are overly complicated. References~\cite{Rauwolf1996, Coverstone2000} both implemented genetic algorithms to solve orbital transfer problems. Reference~\cite{Rauwolf1996} obtained solutions to low-thrust orbit trajectory optimal control problems through using genetic algorithms by testing the effectiveness to determine near-optimal solutions when adding in a thrust/coast arc variable. Reference~\cite{Coverstone2000} obtained solutions to low-thrust orbit trajectory optimal control problems by applying a hybrid optimization method that incorporated a multi-objective genetic algorithm with a low-thrust trajectory optimizer that applies calculus of variations. Reference~\cite{Dachwald2005} combined artificial neural networks and evolutionary algorithms to create evolutionary neurocontrollers in order to solve low-thrust orbital transfer problems. This novel method does not require an initial guess, unlike most traditional optimizers that rely tremendously on an adequate initial guess because of the convergence behavior that causes the optimal solution to be close to the initial guess. References~\cite{Starek2012, Betts2014, Yang2009} implemented model predictive control strategies to solve an orbital transfer problem. Reference~\cite{Starek2012} solved a minimum-fuel Earth-to-Mars rendezvous maneuver to compare two nonlinear model predictive control strategies. This comparison was to showcase the ability to use each strategy to solve problems with interplanetary orbital dynamics and low control authority. The first strategy solved the optimal control problem over a receding horizon with a fixed number of subintervals, which was able to withstand errors in control allocation and unmodeled effects. The second strategy solved the optimal control problem over a receding horizon with a shrinking number of control subintervals, where in order to maintain controllability, a doubling strategy was implemented and shown to handle more complicated bounded problems. Reference~\cite{Betts2014} solved an optimal control problem of a low-thrust orbit transfer with eclipsing. When the spacecraft is passing through the Earth's shadow, the spacecraft is said to be on a coasting arc. A receding horizon algorithm was utilized to create an initial guess and the optimal control problem was divided into separate phases depending on the structure of the thrust. Reference~\cite{Yang2009} created a direct method to find the optimal solution of minimum-time, low-thrust, many-revolution Earth-orbit transfers. A mapping between the Lyapunov control law and the proposed parameterized control law is employed to create a technique that will provide a good initial guess, which had an effect of an enlarged convergence domain of the direct method. References~\cite{JiangBaoyinLi2012, Gergaud2006} both employed a homotopic approach to solve a minimum-fuel low-thrust orbital trajectory optimal control problem. Reference~\cite{JiangBaoyinLi2012} created a method to increase the chance of determining the global optimal solution, while also decreasing the computational time. The first and second-order time derivatives of the switching function were utilized in determining the structure of the control and if it is bang-bang. Reference~\cite{Gergaud2006} implemented a homotopic method to establish convergence properties, where the thrusting structure was found to be discontinuous, which had $1786$ switches that were detected for a thrust of $0.1~\textrm{N}$. Reference~\cite{Daero2021} implemented both indirect and direct methods when solving a minimum-fuel low-thrust Earth-to-Mars orbital transfer. An assumption was made a priori that the thrusting structure would be burn-coast-burn.  For the indirect method, a two-point boundary value problem was derived, whereas for the direct method, an optimal control software $\mathbb{GPOPS-II}$ was used to solve the optimal control problem. References~\cite{KlueverOleson1998, Betts2000, ScheelConway1994, SchubertRao2012, Graham2015} employed a direct optimization method to solve low-thrust orbit transfers. Reference~\cite{KlueverOleson1998} developed a novel direct optimization method for solving near-optimal, minimum-time Earth-orbit transfers when using solar electric propulsion. An optimal blend of three extremal feedback control laws was implemented to create the time history of the optimal thrust direction. This direct method was found to be insensitive to the initial guess and have robust convergence properties. Reference~\cite{Betts2000} explained the implementation of a direct SQP collocation method to the solution of very low-thrust orbit transfers. When using very low-thrust, long duration trajectories are required when the orbital trajectory significantly changes, which results in a very large nonlinear programming problem when using a collocation method. Over 578 revolutions were required for the optimal orbital trajectory, which caused the sparse optimization problem to have 416123 variables and 249674 constraints. Reference~\cite{ScheelConway1994} used a direct transcription approach to solve minimum-time, continuous, very-low-thrust orbit transfers. One of the problems studied was a transfer from low-Earth orbit (LEO) to geostationary orbit (GEO), where modified equinoctial elements were used to describe the motion of the spacecraft. Because of the small thrust accelerations used, the trajectories required multiple revolutions around Earth to attain the terminal conditions of the optimal control problems. Reference~\cite{SchubertRao2012} studied a minimum-time, many-revolutions, low-thrust Earth-orbit trajectory optimization problem using a direct transcription variable-order Gaussian quadrature collocation method. It was discovered that there is a power relationship between the thrust-to-mass ratio and the final mass, transfer time, and final true longitude, where the performance of a given thrust-to-mass ratio can be estimated without solving the optimal control problem. Reference~\cite{Graham2015} solved a variety of minimum-time, low-thrust Earth-orbit transfer problems with the goal of determining high-accuracy solutions. An initial guess was generated by solving a sequence of modified optimal control problems one-revolution at a time until the terminal boundary conditions were satisfied and a variable-order Legendre-Gauss-Radau collocation method was implemented for solving the problems. 

This research is inspired by the work of Ref.~\cite{Herman2002} and both works consider the problem of LEO-to-MEO, LEO-to-HEO, and LEO-to-GEO low-thrust orbital transfer.  It is noted, however, that this research and the work of Ref.~\cite{Herman2002} differ in several significant ways.   First, this research considers seven cases of maximum allowable thrust acceleration for each orbit transfer whereas Ref.~\cite{Herman2002} considers four cases of maximum allowable thrust acceleration.  Second, in this research the thrusting structure is determined as part of the solution process whereas Ref.~\cite{Herman2002} assumed a priori a burn-coast-burn thrusting structure based on the work of Ref.~\cite{Spencer1994}.  In particular, in order to determine the thrusting structure, a recently developed bang-bang and singular optimal control (BBSOC) method \cite{Pager2022} is employed using multi-domain Legendre-Gauss-Radau (LGR) collocation as implemented in the \textsf{MATLAB}$^{\textregistered}$ optimal control software $\mathbb{GPOPS-II}$~\cite{PattersonRao2014}.  It is important to note that, because the thrusting structure is not assumed a priori, the performance (that is, the minimum-impulse) obtained in this study shows a significant improvement over the performance obtained in Ref.~\cite{Herman2002} and this improvement increases as the maximum thrust acceleration decreases.  Third, the improved performance obtained in this study leads to a thrusting structure that is significantly different from the burn-coast-burn thrusting structure assumed in Ref.~\cite{Herman2002}.  Finally, the thrusting structure obtained for the various types of orbital transfer (LEO-to-MEO, LEO-to-HEO, and LEO-to-GEO) is shown for particular cases of maximum allowable thrust acceleration and provide improved insight into the optimal thrusting structure for such transfers.  

This paper is organized as follows. Section~\ref{sect:Problem Description} presents the one-phase minimum-fuel Earth-based orbital transfer optimal control problem by providing the modeling assumptions, dynamics, path constraints, event constraints, the units used to solve the problem, boundary conditions, and variable bounds. Section~\ref{sect:Initial Guess Generation} presents the procedure on how to create the initial guess for each problem case. Section \ref{sect:bbsoc} describes the bang-bang and singular optimal control (BBSOC) methods used to solve the orbital transfer problems in this paper.  Section~\ref{sect:Results and Discussion} shows the key results obtained in this study and provides a discussion of these results.  Finally, Section~\ref{sect:Conclusions} provides conclusions on this research. 

\section{Problem Description}\label{sect:Problem Description}

This section provides the problem description for the orbital transfer problems under consideration.  Section~\ref{sect:Modeling Assumptions} provides all assumptions used to model the motion of the spacecraft.  Section~\ref{sect:Equations of Motion} provides the differential equations of motion for the spacecraft along with the path constraints that are enforced during the motion.  Section~\ref{sect:Scale Factors} provides a description of the scaling that was used to solve the various orbital transfer trajectory optimization problems.  Section~\ref{sect:Boundary Conditions and Bounds} provides the initial and terminal boundary conditions as well as the lower and upper bounds on the time, control, and state.  Finally, Section~\ref{sect:Optimal Control Problem} provides a description of the minimum-fuel orbital transfer optimal control problem.

\subsection{Modeling Assumptions}\label{sect:Modeling Assumptions}

Consider the motion of a point mass in motion over a spherical Earth.  The only forces acting on the spacecraft are those due to thrust and central-body gravitation between the Earth and the spacecraft (all third-body effects and zonal harmonic gravitational coefficients are neglected). Additionally, it is assumed that the control is the thrust magnitude and thrust direction.  

\subsection{Equations of Motion}\label{sect:Equations of Motion}

The dynamics of the spacecraft are modeled using modified equinoctial elements (MEE) $p(t)$, $f(t)$, $g(t)$, $h(t)$, $k(t)$, and $L(t)$ along with the mass, $m(t)$. The relationship between the MEEs and the classical orbital element (COE) are given as follows \cite{MEE}:
\begin{equation}\label{COE-to-MEE}
  \begin{array}{lcl}
    \displaystyle p & = & \displaystyle a(1-e^2), \\
    \displaystyle f & = & \displaystyle e\cos(\omega+\Omega), \\
    \displaystyle g & = & \displaystyle e\sin(\omega+\Omega), \\
    \displaystyle h & = & \displaystyle \tan\left(\frac{i}{2}\right)\cos\Omega, \\
    \displaystyle k & = & \displaystyle \tan\left(\frac{i}{2}\right)\sin\Omega, \\
    \displaystyle L & = & \displaystyle \Omega + \omega + \nu.
  \end{array}
\end{equation}
The differential equations of motion for the spacecraft are given in terms of the MEEs as~\cite{MEE}
\begin{equation}\label{spacecraft-equations-of-motion-wrt-t}
  \begin{array}{lcl}
    \displaystyle \frac{dp}{dt} & = & \displaystyle \sqrt{\frac{p}{\mu_E}}\frac{2p}{w}\Delta_t, \vspace{5pt}\\
    \displaystyle \frac{df}{dt} & = & \displaystyle \sqrt{\frac{p}{\mu_E}}\left[\sin L\Delta_r + \frac{1}{w}\left[\left(w+1\right)\cos L + f\right]\Delta_t - \frac{g}{w}\left(h\sin L - k\cos L\right)\Delta_n\right], \vspace{5pt}\\
    \displaystyle \frac{dg}{dt} & = & \displaystyle \sqrt{\frac{p}{\mu_E}}\left[-\cos L\Delta_r + \frac{1}{w}\left[\left(w+1\right)\sin L + g\right]\Delta_t + \frac{g}{w}\left(h\sin L - k\cos L\right)\Delta_n\right], \vspace{5pt}\\
    \displaystyle \frac{dh}{dt} & = & \displaystyle \sqrt{\frac{p}{\mu_E}}\frac{q^2}{2w}\cos L\Delta_n, \vspace{5pt}\\
    \displaystyle \frac{dk}{dt} & = & \displaystyle \sqrt{\frac{p}{\mu_E}}\frac{q^2}{2w}\sin L\Delta_n, \vspace{5pt}\\
    \displaystyle \frac{dL}{dt} & = & \displaystyle \sqrt{\mu_Ep}\left(\frac{w}{p}\right)^2 + \frac{1}{w}\sqrt{\frac{p}{\mu_E}}\left(h\sin L - k\cos L\right)\Delta_n, \vspace{5pt}\\
    \displaystyle \frac{dm}{dt} & = & \displaystyle -\frac{T}{g_0 I_{sp}}, \vspace{5pt}\\
  \end{array}
\end{equation}
where
\begin{equation}\label{arbitrary-variables}
  \begin{array}{lcl}
    \displaystyle q^2 & = & \displaystyle 1 + h^2 + k^2, \vspace{5pt}\\
    \displaystyle w & = & \displaystyle 1 + f\cos L + g\sin L. \vspace{5pt}\\
  \end{array}
\end{equation}
Due to the thrust of the spacecraft, the non-two-body perturbations in the radial, transverse, and normal directions, $(\Delta_r,\Delta_t,\Delta_n)$, are given as 
\begin{equation}\label{non-two-body-perturbations}
  (\Delta_r,\Delta_t,\Delta_n) = \frac{T}{m} (u_r,u_t,u_n).
\end{equation}
The control consists of the thrust magnitude, $T$, and the thrust direction, $\textbf{u} = (u_r, u_t, u_n)$, where $u_r$ is the component along the direction of the position, $\mathbf{r}$, where $\mathbf{r}$ is measured from the center of the Earth to the spacecraft, $u_n$ is the component along the specific angular momentum, $\mathbf{h}=\mathbf{r}\times\mathbf{v}$ (where $\mathbf{v}$ is the inertial velocity of the spacecraft), and $u_t$ is the component along the direction $\mathbf{h}\times\mathbf{r}$.   Next, in order to ensure that the thrust direction is a unit vector, the control equality path constraint
\begin{equation} \label{thrust-direction-unit-vector}
  \mathbf{u}\cdot\mathbf{u} = u_r^2 + u_t^2 + u_n^2 = 1
\end{equation}
is enforced to guarantee that the thrust direction is a unit vector.  Finally, the following initial and terminal boundary conditions (event constraints) are enforced in order to ensure that the spacecraft starts and terminates in the desired orbits~\cite{Betts2010}:
\begin{equation} \label{event-constraints}
  \begin{array}{lcl}
    \displaystyle p(t_0) & = & a_0 (1-e_0^2), \vspace{5pt}\\
    \displaystyle f^2\left(t_0\right) + g^2\left(t_0\right) & = & \displaystyle e_0^2, \vspace{5pt}\\
    \displaystyle h^2\left(t_0\right) + k^2\left(t_0\right) & = & \displaystyle \tan^2\left(\frac{i_0}{2}\right), \vspace{5pt}\\
    \displaystyle \frac{k\left(t_0\right)}{h\left(t_0\right)} & = & \displaystyle \tan\left(\Omega_0\right), \vspace{5pt}\\
    \displaystyle p(t_f) & = & a_f (1-e_f^2), \vspace{5pt}\\
    \displaystyle f^2\left(t_f\right) + g^2\left(t_f\right) & = & \displaystyle e_f^2, \vspace{5pt}\\
    \displaystyle h^2\left(t_f\right) + k^2\left(t_f\right) & = & \displaystyle \tan^2\left(\frac{i_f}{2}\right), \vspace{5pt}\\
  \end{array}
\end{equation}
It is noted in Eq.~\eqref{event-constraints} that $a_0$, $e_0$, $i_0$, $\Omega_0$, $a_f$, $e_f$, and $i_f$ are parameters determined by the initial and terminal orbits as given in Table~\ref{tab:OrbitalElements}. To model and solve this optimal control problem, the numerical values of the variables required are given in Table~\ref{tab:physical-constants}.
\begin{table}[htbp]
  \centering
  \begin{small}
  \caption{Physical constants.\label{tab:physical-constants}}
  \begin{tabular}{|l|c|} \hline
    Quantity & Value \\ \hline\hline
    $R_E$ & $6.378145\times 10^{6} ~\textrm{m}$\\\hline
    $\mu_E$   & $3.986004418 \times 10^{14} ~\textrm{m}^3 \cdot \textrm{s}^{-2}$   \\\hline
    $g_0$ & $9.80665 \times 10^{5} ~\textrm{m} \cdot \textrm{s}^{-2}$ \\\hline
    $m_0$ & $1.000 \times 10^{3} ~\textrm{kg}$ \\\hline
    $I_{sp}$ & $1.000 \times 10^{3} ~\textrm{s}$ \\\hline
  \end{tabular}
  \end{small}
\end{table}

\subsection{Scale Factors}\label{sect:Scale Factors}

In the optimal control problem, the gravitational parameter of the Earth, $\mu_E$, was set to unity. To attain $\mu_E = 1$, the distance unit (DU), speed unit (VU), time unit (TU), acceleration unit (AU), mass unit (MU), and force unit (FU), were chosen as follows:
\begin{equation}\label{eq:ScaleFactors}
  \left[
    \begin{array}{c}
      \textrm{DU} \\
      \textrm{VU} \\
      \textrm{TU} \\
      \textrm{AU} \\
      \textrm{MU} \\
      \textrm{FU}
    \end{array}
  \right]
  =
  \left[
    \begin{array}{c}
      R_E \\
      \sqrt{\mu_E/\textrm{DU}} \\
      \textrm{DU}/\textrm{VU} \\
      \textrm{VU}/\textrm{TU} \\
      m_0 \\
      \textrm{MU} \cdot \textrm{AU}
    \end{array}
  \right].
\end{equation}

\subsection{Boundary Conditions and Bounds}\label{sect:Boundary Conditions and Bounds}

Table~\ref{tab:OrbitalElements} presents the orbital elements of the initial low-Earth orbit and terminal middle-Earth, high-Earth, or geosynchronous orbits~\cite{Herman2002}. 
\begin{table}[htbp]
 \caption{Orbital elements for the initial and terminal orbits.\label{tab:OrbitalElements}}
 \centering
 \begin{small}
 \begin{tabular}{|l|c|c|c|c|}  \hline
 Orbital Element & LEO & MEO & HEO & GEO \\\hline\hline
 Semi-major Axis, $a~\left(\textrm{km}\right)$ & $7,003$ & $26,560$ & $26,578$ & $42,287$ \\\hline
 Eccentricity, $e$ & $0$ & $0$ & $0.73646$ & $0$ \\\hline
 Inclination, $i~\left(\textrm{deg}\right)$ & $28.5$ & $54.7$ & $63.435$ & $0$ \\\hline
 Longitude of the Ascending Node, $\Omega~\left(\textrm{deg}\right)$ & $0$ & $\textrm{Free}$ & $\textrm{Free}$ & $\textrm{Undefined}$ \\\hline 
 Argument of Periapsis, $\omega~\left(\textrm{deg}\right)$ & $\textrm{Undefined}$ & $\textrm{Undefined}$ & $\textrm{Free}$ & $\textrm{Undefined}$ \\\hline
 True Anomaly, $\nu~\left(\textrm{deg}\right)$ & $\textrm{Free}$ & $\textrm{Free}$ & $\textrm{Free}$ & $\textrm{Free}$\\\hline 
 \end{tabular}
 \end{small}
\end{table}
Bounds are placed on the time, control, and state which are given as
\begin{equation}\label{eq:Bounds}
\begin{array}{lclcll}
t_{0,\min} & \leq & t_0 & \leq & t_{0,\max}\\
t_{f,\min} & \leq & t_f & \leq & t_{f,\max}\\
T_{\min} & \leq & T & \leq & T_{\max} \\
u_{r,\min} & \leq & u_r & \leq & u_{r,\max}\\
u_{t,\min} & \leq & u_t & \leq & u_{t,\max}\\
u_{n,\min} & \leq & u_n & \leq & u_{n,\max}\\
p_{\min} & \leq & p & \leq & p_{\max} \\
f_{\min} & \leq & f & \leq & f_{\max} \\
g_{\min} & \leq & g & \leq & g_{\max} \\
h_{\min} & \leq & h & \leq & h_{\max} \\
k_{\min} & \leq & k & \leq & k_{\max} \\
L_{\min} & \leq & L & \leq & L_{\max} \\
m_{\min} & \leq & m & \leq & m_{\max} \\
\end{array}.\\\\
\end{equation}
All initial and terminal boundary conditions are considered to be free parameters except for the variables specified in Table~\ref{tab:boundary-conditions}. Table~\ref{tab:lower-upper-bounds} shows the variables that are constrained by lower and upper bounds, where all other variables are considered to be free. 
\begin{table}[htbp]
 \caption{Initial and terminal boundary conditions.\label{tab:boundary-conditions}}
 \centering
 \begin{small}
 \begin{tabular}{|l|c|}  \hline
 Variable & Value \\\hline\hline
 $t_0$ & $0~\textrm{d}$ \\\hline
 $p\left(t_0\right)$ & $p_0$ \\\hline
 $m\left(t_0\right)$ & $m_0$\\\hline
 $p\left(t_f\right)$ & $p_f$ \\\hline
 \end{tabular}
 \end{small}
\end{table}
\begin{table}[htbp]
 \caption{Lower and upper bounds.\label{tab:lower-upper-bounds}}
 \centering
 \begin{small}
 \begin{tabular}{|l|c|}  \hline
 Variable & [Lower, Upper]\\\hline\hline
 $t$ & $[0, \textrm{Free}]$ \\\hline
 $T$ & $[0, T_{\max}]$ \\\hline
 $f$ & $[-1, +1]$ \\\hline
 $g$ & $[-1, +1]$ \\\hline
 $h$ & $[-1, +1]$ \\\hline
 $k$ & $[-1, +1]$ \\\hline
 $m$ & $[0.01*m_0, m_0]$ \\\hline
 \end{tabular}
 \end{small}
\end{table}
It is noted that $p_0$ and $p_f$ correspond to the initial and terminal semi-parameter, respectively, which can be calculated from the values in Table~\ref{tab:OrbitalElements} and the relationship $p=a(1-e^2)$ from Eq.~\ref{COE-to-MEE}. Finally, Table~\ref{tab:thrust-values} gives the maximum allowable thrust acceleration, $s_0$, and the maximum thrust magnitude, $T_{\max}$, for all cases studied in this research~\cite{Herman2002}.
\begin{table}[htbp]
  \caption{Maximum allowable thrust acceleration and maximum thrust values for all cases.\label{tab:thrust-values}}
  \centering
  \begin{small}
  \begin{tabular}{|c|c|c|} \hline
    Case & $s_0\left(\textrm{m} \cdot \textrm{s}^{-2}\right)$ & $T_{\max}\left(\textrm{N}\right)$ \\\hline
    1 & $10$    & $10000$ \\\hline
    2 & $5$      & $5000$ \\\hline
    3 & $1$      & $1000$ \\\hline
    4 & $0.5$   & $500$ \\\hline
    5 & $0.1$   & $100$ \\\hline
    6 & $0.05$ & $50$ \\\hline
    7 & $0.01$ & $10$ \\\hline
  \end{tabular}
  \end{small}
\end{table}

\subsection{Optimal Control Problem}\label{sect:Optimal Control Problem}

The optimal control problem for the various Earth orbit transfers is stated as follows. Determine the state $\left(p,f,g,h,k,L,m\right)$, control $\left(T,u_r,u_t,u_n\right)$, and terminal time $t_f$ which maximizes the final mass of the spacecraft, $m\left(t_f\right)$. For that reason, the following cost functional needs to be minimized
\begin{equation}\label{eq:objective-functional}
J = -m\left(t_f\right)
\end{equation}
while satisfying the dynamics and constraints in Section~\ref{sect:Equations of Motion} along with the variable bounds and boundary conditions in Section~\ref{sect:Boundary Conditions and Bounds}.

\section{Initial Guess Generation}\label{sect:Initial Guess Generation}

This section explains how initial guesses required by the general-purpose \textsf{MATLAB}${}^{\textregistered}$ optimal control software $\mathbb{GPOPS-II}$~\cite{PattersonRao2014} are generated to solve the single-phase minimum-fuel Earth-based orbit transfer optimal control problem described in Section~\ref{sect:Problem Description}. Two distinct initial guess generation methods are described in Section~\ref{sect:Initial Guess for Partial Orbital Revolution Solutions} and Section~\ref{sect:Initial Guess for Multiple Orbital Revolution Solutions}.  The first initial guess generation method is used for maximum allowable thrust acceleration values that lead to solutions that are less than one orbital revolution (that is, a partial orbital revolution) while the second initial guess generation method is used for maximum allowable thrust acceleration values that lead to solutions that are more than one orbital revolution (that is, multiple orbital revolutions).   In order to determine which generation method is applied, each case of the optimal control problem is solved on the $\mathbb{GPOPS-II}$ default initial mesh (that is, no mesh refinement is implemented) with the initial guess generation method described in Section~\ref{sect:Initial Guess for Partial Orbital Revolution Solutions} to determine if the optimal solution is categorized as either a partial orbital revolution or a multiple orbital revolution solution. If the solution is a partial orbital revolution, then the optimal control problem is solved using the guess generation method in Section~\ref{sect:Initial Guess for Partial Orbital Revolution Solutions}.  On the other hand, if the solution consists of multiple orbital revolutions, then the optimal control problem is solved using the generation method in Section~\ref{sect:Initial Guess for Multiple Orbital Revolution Solutions}.

\subsection{Initial Guess for Partial Orbital Revolution Solutions}\label{sect:Initial Guess for Partial Orbital Revolution Solutions}

For the problem cases categorized as partial orbital revolution solutions, the ordinary differential equations solver \textsf{ode113} in \textsf{MATLAB}${}^{\textregistered}$ is used to generate the initial guess. The ODE solver integrates the spacecraft dynamics, $(p\left(L\right),f\left(L\right),g\left(L\right),h\left(L\right),k\left(L\right),t\left(L\right),m\left(L\right))$, in Eqs.~\eqref{t-to-L-conversion} and~\eqref{spacecraft-equations-of-motion-wrt-L}. In this procedure, the true longitude, $L$, replaces time, $t$, as the independent variable, therefore the dynamics in Eq.~\eqref{spacecraft-equations-of-motion-wrt-t} must be transformed using the conversion factor
\begin{equation} \label{t-to-L-conversion}
\begin{array}{lcl}
  \displaystyle \frac{dt}{dL} & = & \displaystyle \left(\frac{dL}{dt}\right)^{-1}, \vspace{5pt}\\
\end{array}
\end{equation}
so that the other states of the spacecraft, $(p(L),f(L),g(L),h(L),k(L),m(L))$, are given as 
\begin{equation} \label{spacecraft-equations-of-motion-wrt-L}
\begin{array}{lcl}
  \displaystyle \frac{dp}{dL} & = & \displaystyle \frac{dt}{dL}\frac{dp}{dt}, \vspace{5pt}\\
  \displaystyle \frac{df}{dL} & = & \displaystyle \frac{dt}{dL}\frac{df}{dt}, \vspace{5pt}\\
  \displaystyle \frac{dg}{dL} & = & \displaystyle \frac{dt}{dL}\frac{dg}{dt}, \vspace{5pt}\\
  \displaystyle \frac{dh}{dL} & = & \displaystyle \frac{dt}{dL}\frac{dh}{dt}, \vspace{5pt}\\
  \displaystyle \frac{dk}{dL} & = & \displaystyle \frac{dt}{dL}\frac{dk}{dt}, \vspace{5pt}\\
  \displaystyle \frac{dm}{dL} & = & \displaystyle \frac{dt}{dL}\frac{dm}{dt}. \vspace{5pt}\\
\end{array}
\end{equation}
The initial conditions are set as the LEO orbital elements in Table~\ref{tab:OrbitalElements} that have been converted to modified equinoctial elements. The initial guess is integrated until the semi-parameter, $p$, of the orbit is equal to the corresponding terminal value, $p_f$. For the control, $\left(T,u_r,u_t,u_n\right)$, the thrust magnitude, $T$, is considered to be at maximum for the entirety of the integration and the thrust direction components are solved for afterwards by assuming that the thrust is always in the same direction as the velocity vector.

\subsection{Initial Guess for Multiple Orbital Revolution Solutions}\label{sect:Initial Guess for Multiple Orbital Revolution Solutions}

For the problem cases categorized as multiple orbital revolution solutions, an initial guess procedure was implemented~\cite{Graham2015}. The cases that contain multiple orbital revolutions require an initial guess that contains a number of orbital revolutions that is reasonably close to the actual number of orbital revolutions of the optimized solution that the NLP solver converged to. This initial guess procedure consists of solving a chain of optimal control sub-problems, where the goal is to minimize the following objective functional that consists of a mean square relative difference
\begin{equation}\label{Initial-Guess-Objective-Functional}
J = \left[\frac{p\left(L_f\right)-p_D}{1 + p_D}\right]^2 + \left[\frac{f^2\left(L_f\right)+g^2\left(L_f\right)-e_D^2}{1 + e_D^2}\right]^2 + \left[\frac{h^2\left(L_f\right)+k^2\left(L_f\right)-\tan^2\left(\frac{i_D}{2}\right)}{\sec^2\left(\frac{i_D}{2}\right)}\right]^2,
\end{equation}
by determining the state and control that transfers the spacecraft from the initial low-Earth orbit to the correct terminal orbit, depending on the study. When minimizing Eq.~\eqref{Initial-Guess-Objective-Functional}, the sub-problem attains an optimal solution that is as close in proximity as possible to the desired terminal semi-parameter, $p_D$, eccentricity, $e_D$, and inclination, $i_D$, within one orbital revolution. Each sub-problem uses the terminal state of the previous sub-problem as the initial state of the current sub-problem and is evaluated at most over one orbital revolution. It is noted that for the first sub-problem, the initial state is set as the LEO orbital elements in Table~\ref{tab:OrbitalElements} that have been converted into modified equinoctial elements using Eq.~\eqref{COE-to-MEE}. In this procedure, the true longitude, $L$, replaces time, $t$, as the independent variable, therefore the dynamics in Eq.~\eqref{spacecraft-equations-of-motion-wrt-t} must be transformed using the conversion factor in Eq.~\eqref{t-to-L-conversion}, so that the state of the spacecraft, $(p(L),f(L),g(L),h(L),k(L),t(L),m(L))$, is given in Eqs.~\eqref{t-to-L-conversion} and~\eqref{spacecraft-equations-of-motion-wrt-L}. 

The continuous-time optimal control sub-problem is then stated as follows. Minimize the objective functional in Eq.~\eqref{Initial-Guess-Objective-Functional}, subject to the dynamics constraints in Eqs.~\eqref{t-to-L-conversion} and~\eqref{spacecraft-equations-of-motion-wrt-L}, the path constraint in Eq.~\eqref{thrust-direction-unit-vector}, and the boundary conditions 
\begin{equation}\label{Initial-Guess-Boundary-Conditions}
\begin{array}{rclcrcl}\vspace{5pt}
p^{\left(n\right)}\left(L_0^{\left(n\right)}\right) & = & p^{\left(n-1\right)}\left(L_f^{\left(n-1\right)}\right), & & p^{\left(n\right)}\left(L_f^{\left(n\right)}\right) & = & \textrm{Free}, \\\vspace{5pt}

f^{\left(n\right)}\left(L_0^{\left(n\right)}\right) & = & f^{\left(n-1\right)}\left(L_f^{\left(n-1\right)}\right), & & f^{\left(n\right)}\left(L_f^{\left(n\right)}\right) & = & \textrm{Free}, \\\vspace{5pt}

g^{\left(n\right)}\left(L_0^{\left(n\right)}\right) & = & g^{\left(n-1\right)}\left(L_f^{\left(n-1\right)}\right), & & g^{\left(n\right)}\left(L_f^{\left(n\right)}\right) & = & \textrm{Free}, \\\vspace{5pt}

h^{\left(n\right)}\left(L_0^{\left(n\right)}\right) & = & h^{\left(n-1\right)}\left(L_f^{\left(n-1\right)}\right), & & h^{\left(n\right)}\left(L_f^{\left(n\right)}\right) & = & \textrm{Free}, \\\vspace{5pt}

k^{\left(n\right)}\left(L_0^{\left(n\right)}\right) & = & k^{\left(n-1\right)}\left(L_f^{\left(n-1\right)}\right), & & k^{\left(n\right)}\left(L_f^{\left(n\right)}\right) & = & \textrm{Free}, \\\vspace{5pt}

t^{\left(n\right)}\left(L_0^{\left(n\right)}\right) & = & t^{\left(n-1\right)}\left(L_f^{\left(n-1\right)}\right), & & t^{\left(n\right)}\left(L_f^{\left(n\right)}\right) & = & \textrm{Free}, \\\vspace{5pt}

m^{\left(n\right)}\left(L_0^{\left(n\right)}\right) & = & m^{\left(n-1\right)}\left(L_f^{\left(n-1\right)}\right), & & m^{\left(n\right)}\left(L_f^{\left(n\right)}\right) & = & \textrm{Free}, \\\vspace{5pt}

L_0^{\left(n\right)} & = & L_f^{\left(n-1\right)}, & & L_f^{\left(n\right)} & \leq & L_f^{\left(n-1\right)} + 2\pi, \\\vspace{5pt}
\end{array}
\end{equation}
for $n = 1,\cdots,N$, where $N$ is the total number of true longitude cycles. The sub-problem solutions are combined into an initial guess once the desired terminal conditions of the terminal orbital elements in Table~\ref{tab:OrbitalElements} are within an error tolerance of $10^{-4}$.

The sub-problems are solved using the general-purpose \textsf{MATLAB}${}^{\textregistered}$ optimal control software $\mathbb{GPOPS-II}$~\cite{PattersonRao2014} together with the nonlinear program (NLP) solver {\em IPOPT}~\cite{BieglerZavala2009} in full Newton (second derivative) mode with an NLP solver tolerance of $10^{-7}$.  All derivatives required by {\em IPOPT} were obtained using the open-source algorithmic differentiation software {\em ADiGator}~\cite{WeinsteinRao2017}.   $\mathbb{GPOPS-II}$ was employed using the following settings.  First, an initial guess is created using the ordinary differential equations solver \textsf{ode113}, where the ODE solver integrates the spacecraft dynamics, $(p\left(L\right),f\left(L\right),g\left(L\right),h\left(L\right),k\left(L\right),t\left(L\right),m\left(L\right))$, in Eqs.~\eqref{t-to-L-conversion} and~\eqref{spacecraft-equations-of-motion-wrt-L}. The initial conditions are set as the LEO orbital elements in Table~\ref{tab:OrbitalElements} that have been converted to modified equinoctial elements using Eq.~\eqref{COE-to-MEE} for the first sub-problem and then, for every subsequent sub-problem, the initial conditions are set as the terminal values of the previous sub-problem. The initial guess is integrated for one orbital revolution, where the thrust magnitude, $T$, is considered to be at maximum for the entirety of the integration and the thrust direction components, $\left(u_r,u_t,u_n\right)$, are solved for afterwards by assuming that the thrust is always in the same direction as the velocity vector. Second, a $hp$-adaptive mesh refinement method~\cite{LiuHagerRao2018} is used with a mesh refinement accuracy tolerance of $10^{-2}$. Third, the initial mesh is set to have one mesh interval with four collocation points.  Finally, all computations were performed using a 2.9 GHz 6-Core Intel Core i9 MacBook Pro running Mac OS version 11.6.1 (Big Sur) with 32GB 2400MHz DDR4 RAM and \textsf{MATLAB}${}^{\textregistered}$ Version R2018b (build 9.5.0.944444).

\section{Numerical Approach: BBSOC Method}\label{sect:bbsoc}

The minimum-fuel orbital transfer problems described in Section \ref{sect:Problem Description} are solved using a recently developed bang-bang and singular optimal control (BBSOC) method developed in Ref.~\cite{Pager2022}. The BBSOC method is capable of identifying the existence of both bang-bang switches in the optimal control along with identifying singular arcs.  For segments identified as singular arcs, the BBSOC method performs an iterative regularization procedure to compute the singular control. For segments identified as bang-bang, the BBSOC method determines whether the control lies at either its lower or upper limit and optimize values of the switch times. Specifically, the structure identification procedure of the BBSOC method analyzes an initial solution to an optimal control problem for any discontinuities, bang-bang arcs, and/or singular arcs. After the structure has been identified, the BBSOC method partitions the initial mesh into domains representing each identified arc.  Finally, corresponding constraints are enforced in each domain, respectively, while a regularization procedure is applied in any domains that have been identified as singular. The method is algorithmic in nature and requires no user-input.  In addition, the BBSOC method is implemented using a multiple-domain formulation of the $hp$-adaptive Legendre-Gauss-Radau (LGR) collocation method~\cite{Garg2010,Garg2011a,Garg2011b,Kameswaran2008,PattersonHagerRao2015}, which allows the use of variable mesh refinement to exploit knowledge of the identified control structure. More details on the BBSOC method and the multiple-domain LGR collocation can be found in Refs.~\cite{Pager2022}.

\section{Results and Discussion}\label{sect:Results and Discussion}

This section discusses the results obtained by solving each Earth-based orbital transfer study, where the optimal control problem is described in Section~\ref{sect:Problem Description}. Each of the three studies are solved for seven cases of maximum allowable thrust acceleration, $s_0$, with the values shown in Table~\ref{tab:thrust-values}.  Section~\ref{sect:overall-performance} provides an analysis of the overall performance for all three types of orbital transfer, while Section~\ref{sect:key-features} discusses the key features of the optimized solutions for each type of orbital transfer.
All results presented in this section were obtained using the bang-bang and singular optimal control (BBSOC) method~\cite{Pager2022} implemented in \textsf{MATLAB}\textsuperscript{\textregistered} with the NLP problem solver IPOPT~\cite{BieglerZavala2009} employed in full Newton (second derivative) mode with an NLP solver tolerance of $10^{-7}$.  In addition, the method of Ref.~\cite{PattersonHagerRao2015} is used for any necessary mesh refinement.  All derivatives required by IPOPT are computed using the algorithmic differentiation software {\em ADiGator}~\cite{WeinsteinRao2017}. Finally, all computations were performed using a 2.9 GHz 6-Core Intel Core i9 MacBook Pro running Mac OS version 11.6.1 (Big Sur) with 32GB 2400MHz DDR4 RAM and \textsf{MATLAB}\textsuperscript{\textregistered} Version R2018b (build 9.5.0.944444).

Tables~\ref{tab:Study 1: InitialSetup} -~\ref{tab:Study 3: InitialSetup} show the initial setup for the BBSOC method for solving the optimal control problems of Study 1: LEO-to-MEO; Study 2: LEO-to-HEO; and Study 3: LEO-to-GEO, respectively, where $\eta$ is a threshold to detect the relative size of the jumps in the control, and $M$ is the number of initial mesh intervals. It is noted that, for every case, the number of initial collocation points in each initial mesh interval is $c_p = 3$.
\begin{table}[htbp]
 \caption{Initial setup for each case of Study 1: LEO-to-MEO.\label{tab:Study 1: InitialSetup}}
 \centering
 \begin{small}
 \begin{tabular}{|c|c|c|}  \hline
 Case & $\eta$ & $M$ \\\hline\hline
 1 & 0.1   & 130 \\\hline
 2 & 0.1   & 100 \\\hline
 3 & 0.1   & 10   \\\hline
 4 & 0.01 & 60   \\\hline
 5 & 0.01 & 50   \\\hline
 6 & 0.01 & 40   \\\hline
 7 & 0.1   & 50   \\\hline
 \end{tabular}
 \end{small}
\end{table}
\begin{table}[htbp]
 \caption{Initial setup for each case of Study 2: LEO-to-HEO.\label{tab:Study 2: InitialSetup}}
 \centering
 \begin{small}
 \begin{tabular}{|c|c|c|}  \hline
 Case & $\eta$ & $M$ \\\hline\hline
 1 & 0.01   & 170 \\\hline
 2 & 0.001 & 70   \\\hline
 3 & 0.01   & 90   \\\hline
 4 & 0.01   & 50   \\\hline
 5 & 0.01   & 30   \\\hline
 6 & 0.001 & 140 \\\hline
 7 & 0.1     & 60   \\\hline
 \end{tabular}
 \end{small}
\end{table}
\begin{table}[htbp]
 \caption{Initial setup for each case of Study 3: LEO-to-GEO.\label{tab:Study 3: InitialSetup}}
 \centering
 \begin{small}
 \begin{tabular}{|c|c|c|}  \hline
 Case & $\eta$ & $M$ \\\hline\hline
 1 & 0.1     & 90   \\\hline
 2 & 0.1     & 190 \\\hline
 3 & 0.01   & 10   \\\hline
 4 & 0.01   & 50   \\\hline
 5 & 0.01   & 40   \\\hline
 6 & 0.001 & 150 \\\hline
 7 & 0.1     & 110 \\\hline
 \end{tabular}
 \end{small}
\end{table}

\subsection{Overall Performance of Optimized Orbital Transfers\label{sect:overall-performance}}

Tables~\ref{tab:Study 1: Overall Results}--\ref{tab:Study 3: Overall Results} show the overall performance for the LEO-to-MEO, LEO-to-HEO, and LEO-to-GEO transfers, respectively.  In particular, Tables~\ref{tab:Study 1: Overall Results}--\ref{tab:Study 3: Overall Results} show the terminal mass, $m(t_f)$, the total thrusting time, $t_{T}$, the total number of revolutions, $N$, the number of thrust arcs, $A_{T}$, and the total impulse, $\Delta V$, for the seven values of maximum allowable thrust acceleration, $s_0$, studied. In addition, for those values of $s_0$ that appear in Ref.~\cite{Herman2002}, the impulse, $\Delta V$, obtained in Ref.~\cite{Herman2002} is shown in Tables~\ref{tab:Study 1: Overall Results}--\ref{tab:Study 3: Overall Results}. When computing the impulse, $\Delta V$, the Tsiolkovsky rocket equation
\begin{equation}\label{RocketEquation}
\Delta V = g_0 I_{sp} \ln \left[ \frac{m_0}{m\left(t_f\right)} \right].
\end{equation}
is used. For all three types of transfers the solutions can be separated into two categories: orbital transfers with partial orbital revolutions and orbital transfers with multiple orbital revolutions. In particular, regardless of the type of transfer, the partial revolution transfers correspond to larger values of $s_0$ (Cases 1-4), whereas the multiple revolution transfers correspond to smaller values of $s_0$ (Cases 5-7). Furthermore, as $s_0$ decreases, the terminal mass decreases, the total time thrusting increases, the total revolutions increase, and the total impulse increases.  It is noted that, for the partial orbital revolution solutions, each case has two thrust arcs whereas, for the multiple orbital revolution solutions, the number of thrust arcs increases as $s_0$ decreases. It is noted that the terminal mass decreases as $s_0$ decreases because the total time thrusting increases which leads to an increase in fuel consumption in order to attain the required terminal conditions. When comparing the total impulse of this work to the total impulse obtained in Ref.~\cite{Herman2002} (that is, for the four overlapping cases of maximum allowable thrust acceleration), it can be seen that this work produces a smaller $\Delta V$ (except for Case 1 of Study 2). For Cases 1 and 3, the differences in $\Delta V$ are relatively small because the solutions for these cases consist of only two thrust arcs (that is, burn-coast-burn) which is similar to the burn-coast-burn thrusting structure assumed in Ref.~\cite{Herman2002}. Note, however, that for Cases 1 and 3 the burn-coast-burn thrusting structure is obtained algorithmically via the BBSOC method without any a priori assumptions. On the other hand, for Cases 5 and 7 the differences in $\Delta V$ are significantly larger because the optimal control consists of significantly more than two thrust arcs which thereby shows that the burn-coast-burn thrusting structure assumed in Ref.~\cite{Herman2002} is less fuel efficient. Therefore, the $\Delta V$ obtained in this research is similar to the $\Delta V$ obtained in Ref.~\cite{Herman2002} when the BBSOC method produces a thrusting structure similar to the assumed thrusting structure in Ref.~\cite{Herman2002} and produces a significantly lower $\Delta V$ as the BBSOC method produces a thrusting structure that significantly differs from the assumed thrusting structure of Ref.~\cite{Herman2002}.

{\noindent \bf \em Aside:} It is noted for the LEO-to-HEO transfers that the solutions obtained for Cases 1 and 2 require that a procedure different from that described in Section~\ref{sect:Initial Guess Generation} is used. First, a solution is obtained for Case 2 using the initial setup $(\eta,M,c_p)=(0.01,20,3)$.  As it turns out, the solution obtained on the first mesh using this setup satisfies the accuracy tolerance.  Consequently, the BBSOC method never executes structure detection. In order to make it possible for the BBSOC method to execute structure detection, a modified initial setup $\eta = 0.1$, $M = 10$ and $c_p = 3$ is used. The solution obtained using this second setup is then used as an initial guess with yet a third initial setup $(\eta,M,c_p)=(0.001,70,3)$, leading to the results found in Table~\ref{tab:Study 2: Overall Results}. The solution for Case 2 shown in Table~\ref{tab:Study 2: Overall Results} is then used as an initial guess for Case 1 using the initial setup $(\eta,M,c_p)=(0.1,10,3)$.   The solution obtained using this last setup is then used as an initial guess with the setup $(\eta,M,c_p)=(0.01,170,3)$, leading to the results found in Table~\ref{tab:Study 2: Overall Results}. 

\begin{table}[htbp]
  \centering
  \caption{Performance results for Study 1: LEO-to-MEO transfers.\label{tab:Study 1: Overall Results}}
  \begin{small}
  \begin{tabular}{|c|c|c|c|c|c|c|c|} \hline
    Case & $s_0~\left(\textrm{m} \cdot \textrm{s}^{-2}\right)$ & $m\left(t_f\right)~(\textrm{kg})$ & $t_{T}~\left(\textrm{h}\right)$ & $N$ & $A_{T}$ & $\Delta V~\left(\textrm{m} \cdot \textrm{s}^{-1}\right)$ & Ref.~\cite{Herman2002} $\Delta V~\left(\textrm{m} \cdot \textrm{s}^{-1}\right)$ \\\hline\hline 
    1 & $10$   & $674.9651$ & $0.0880$     & $0.5196$   & $2$ & $3854.9$ & $3863$ \\\hline
    2 & $5$   & $674.7867$ & $0.1744$     & $0.5352$   & $2$ & $3857.5$ & -- \\\hline    
    3 & $1$   & $668.2949$ & $0.8941$     & $0.7398$   & $2$ & $3952.3$ & $3970$ \\\hline    
    4 & $0.5$ & $653.0154$ & $1.8816$     & $0.8768$   & $2$ & $4179.1$ & -- \\\hline
    5 & $0.1$ & $624.2352$ & $10.2132$   & $4.9579$   & $4$ & $4621.2$ & $4731$ \\\hline
    6 & $0.05$ & $607.2275$ & $21.2975$   & $9.1414$   & $6$ & $4892.1$ & -- \\\hline
    7 & $0.01$ & $606.9697$ & $106.4970$ & $32.8742$ & $9$ & $4896.2$ & $5122$ \\\hline
  \end{tabular}
  \end{small}
\end{table}

\begin{table}[htbp]
  \centering
  \caption{Performance results for Study 2: LEO-to-HEO transfers.\label{tab:Study 2: Overall Results}}
  \begin{small}
  \begin{tabular}{|c|c|c|c|c|c|c|c|} \hline
    Case & $s_0~\left(\textrm{m} \cdot \textrm{s}^{-2}\right)$ & $m\left(t_f\right)~(\textrm{kg})$ & $t_{T}~\left(\textrm{h}\right)$ & $N$ & $A_{T}$ & $\Delta V~\left(\textrm{m} \cdot \textrm{s}^{-1}\right)$ & Ref.~\cite{Herman2002} $\Delta V~\left(\textrm{m} \cdot \textrm{s}^{-1}\right)$ \\\hline\hline 
    1 & $10$   & $716.2925$ & $0.0766$     & $0.5202$   & $2$   & $3272.2$ &  $3271$ \\\hline
    2 & $5$   & $715.7138$ & $0.1525$     & $0.5387$   & $2$   & $3280.1$ &  --         \\\hline    
    3 & $1$   & $699.2824$ & $0.8140$     & $0.8359$   & $2$   & $3507.8$ &  $3555$ \\\hline    
    4 & $0.5$ & $663.1665$ & $1.8242$     & $0.9240$   & $2$   & $4027.9$ &  --          \\\hline
    5 & $0.1$ & $657.2695$ & $9.2494$     & $4.9570$   & $6$   & $4115.6$ &  $5271$ \\\hline
    6 & $0.05$ & $645.0881$ & $19.1564$   & $9.0376$   & $9$   & $4298.9$ &  --          \\\hline
    7 & $0.01$ & $576.7825$ & $114.7384$ & $39.0413$ & $17$ & $5396.5$ &  $6109$ \\\hline
  \end{tabular}
  \end{small}
\end{table}

\begin{table}[htbp]
  \centering
  \caption{Performance results for Study 3: LEO-to-GEO transfers.\label{tab:Study 3: Overall Results}}
  \begin{small}
  \begin{tabular}{|c|c|c|c|c|c|c|c|} \hline
    Case & $s_0~\left(\textrm{m} \cdot \textrm{s}^{-2}\right)$ & $m\left(t_f\right)~(\textrm{kg})$ & $t_{T}~\left(\textrm{h}\right)$ & $N$ & $A_{T}$ & $\Delta V~\left(\textrm{m} \cdot \textrm{s}^{-1}\right)$ & Ref.~\cite{Herman2002} $\Delta V~\left(\textrm{m} \cdot \textrm{s}^{-1}\right)$ \\\hline\hline 
    1 & $10$   & $656.7935$ & $0.0925$     & $0.5195$     & $2$ & $4122.6$ & $4127$ \\\hline
    2 & $5$   & $656.3850$ & $0.1857$     & $0.5408$     & $2$ & $4128.7$ & -- \\\hline    
    3 & $1$   & $646.4416$ & $0.9614$     & $0.7694$     & $2$ & $4278.4$ & $4308$ \\\hline    
    4 & $0.5$ & $626.2787$ & $2.0140$     & $0.9380$     & $2$ & $4589.1$ & -- \\\hline
    5 & $0.1$ & $619.0090$ & $10.3158$   & $4.8044$     & $5$ & $4703.6$ & $5167$ \\\hline
    6 & $0.05$ & $583.7997$ & $22.5498$   & $8.0286$     & $6$ & $5277.9$ & -- \\\hline
    7 & $0.01$ & $579.8979$ & $114.0104$ & $110.0091$ & $8$ & $5343.7$ & $5698$ \\\hline
  \end{tabular}
  \end{small}
\end{table}

\subsection{Key Features of Optimized Solutions\label{sect:key-features}}

In this section the key features of the optimized LEO-to-MEO, LEO-to-HEO, and LEO-to-GEO transfers is analyzed by studying the solution for Case 5, $s_0 = 0.1~\textrm{m} \cdot \textrm{s}^{-2}$.  Each type of transfer is studied separately in order to highlight the key features of the solution for that type of transfer.  

\subsubsection{Key Features of Optimized LEO-to-MEO Transfers}\label{sect: Study 1: LEO-to-MEO Key Features of Optimized Solutions}

\begin{figure}[htb]
\centering
\includegraphics[width=3in]{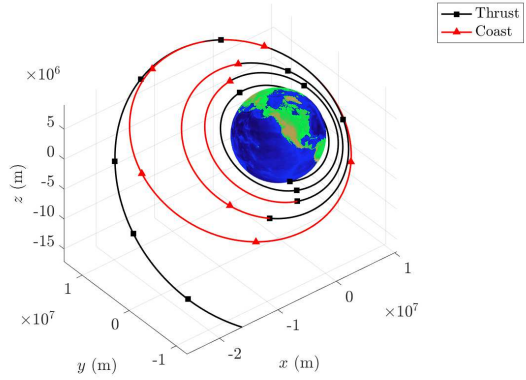}
\caption{Optimal three-dimensional spacecraft trajectory for LEO-to-MEO transfer with $s_0 = 0.1~\textrm{m} \cdot \textrm{s}^{-2}$.\label{fig:MEO_Tm-1_3D_Trajectory}}
\end{figure}

Figure~\ref{fig:MEO_Tm-1_3D_Trajectory} shows the optimized three-dimensional trajectory of the LEO-to-MEO transfer, where the modified equinoctial elements were converted into scaled Cartesian coordinates~\cite{MEE}. It is seen that the spacecraft begins in a low-Earth orbit and terminates in a middle-Earth orbit that corresponds to the orbital elements in Table~\ref{tab:OrbitalElements}. The optimal trajectory of the spacecraft consists of $4.9579$ orbital revolutions around the Earth with a final mass of $624.2352~\textrm{kg}$, four thrust arcs, a total time thrusting of $10.2132~\textrm{h}$, and a total impulse of $4621.2~\textrm{m} \cdot \textrm{s}^{-1}$. It is noted that all of the thrust arcs, except for the final, occur near the periapsis of the orbit. 

\begin{figure}[h]
  \centering
  \subfloat[Semi-major axis.\label{fig:MEO_Tm-1_avst}]{\includegraphics[width=3in]{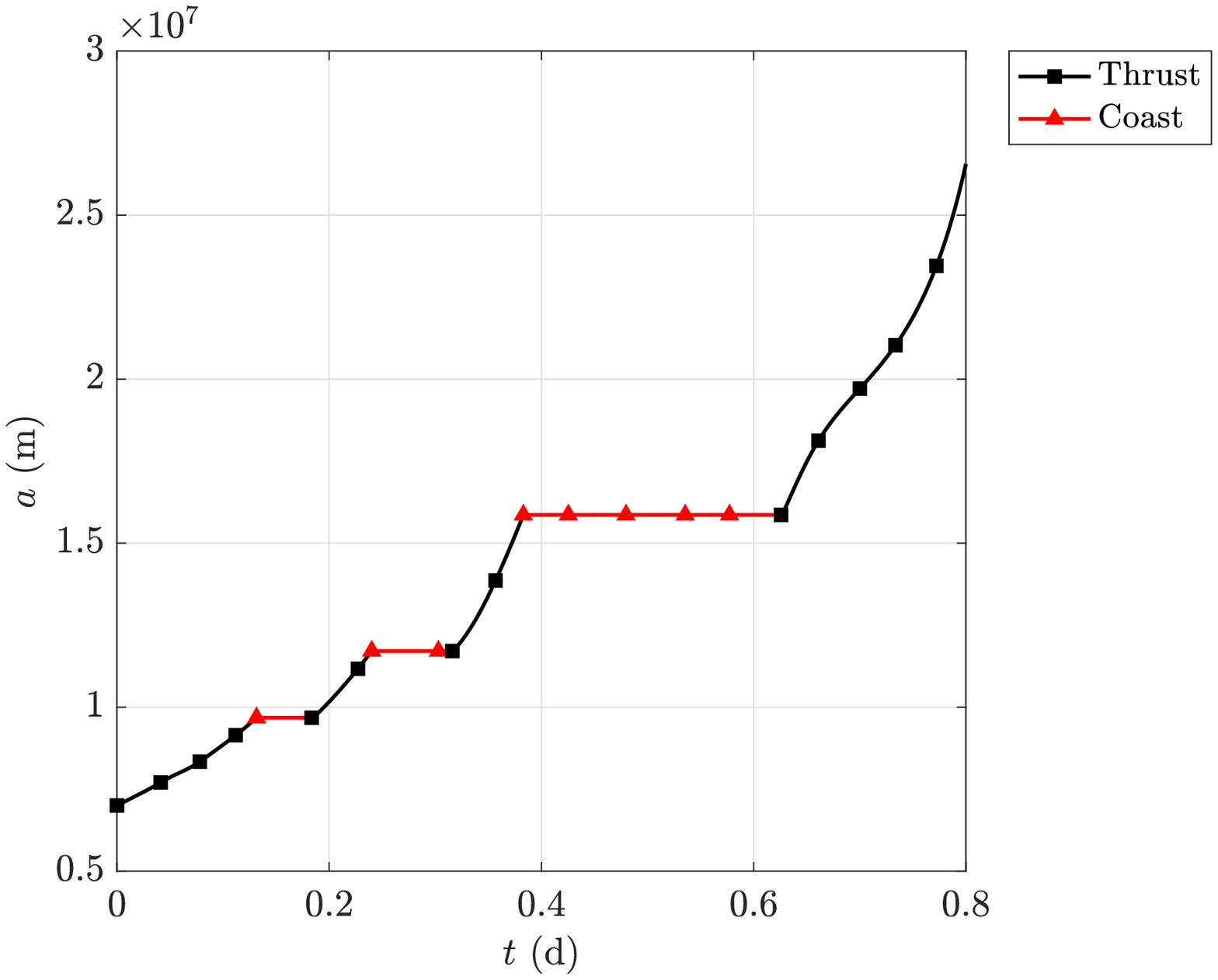}}~~\subfloat[Eccentricity. \label{fig:MEO_Tm-1_evst}]{\includegraphics[width=3in]{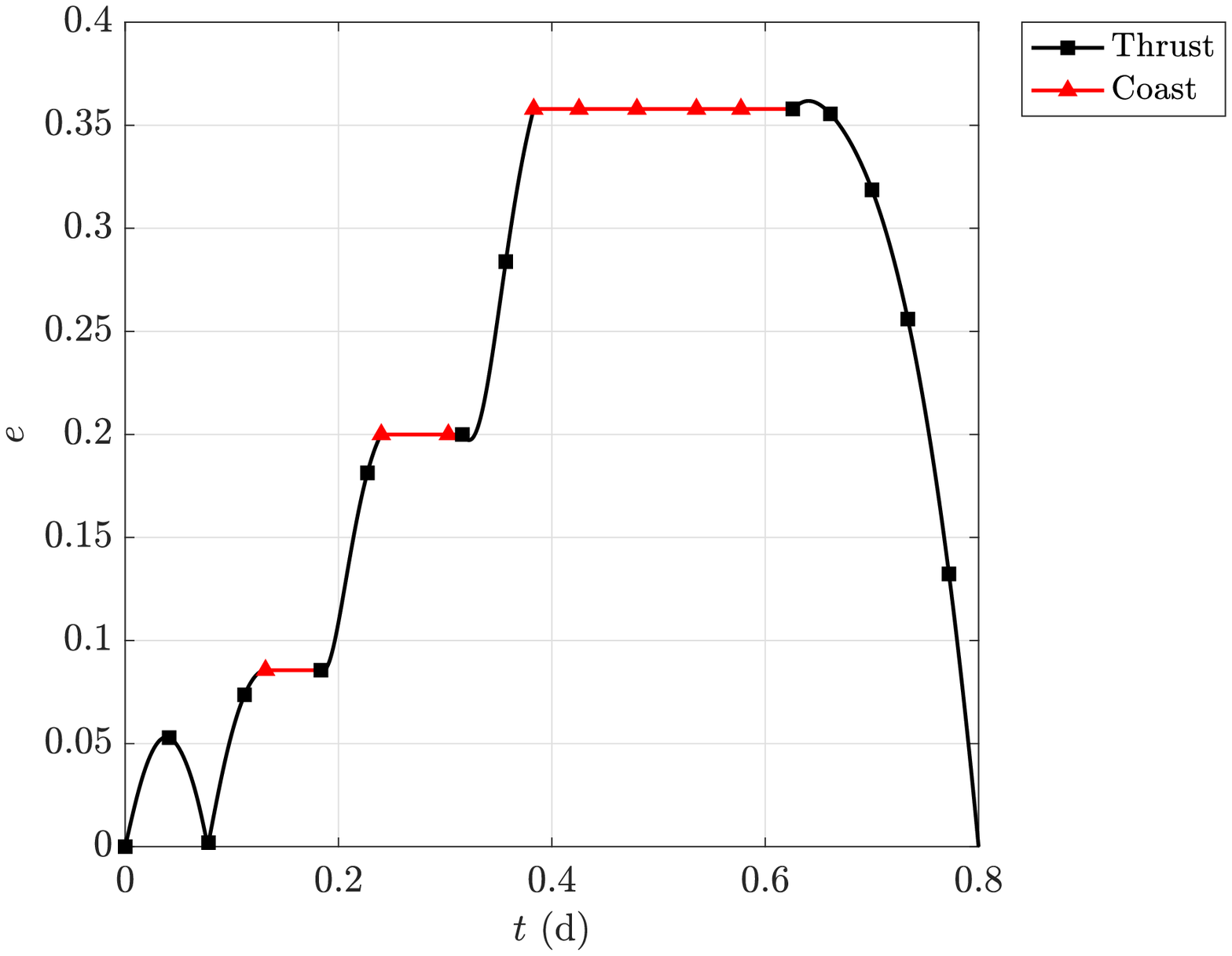}}

  \subfloat[Inclination. \label{fig:MEO_Tm-1_ivst}]{\includegraphics[width=3in]{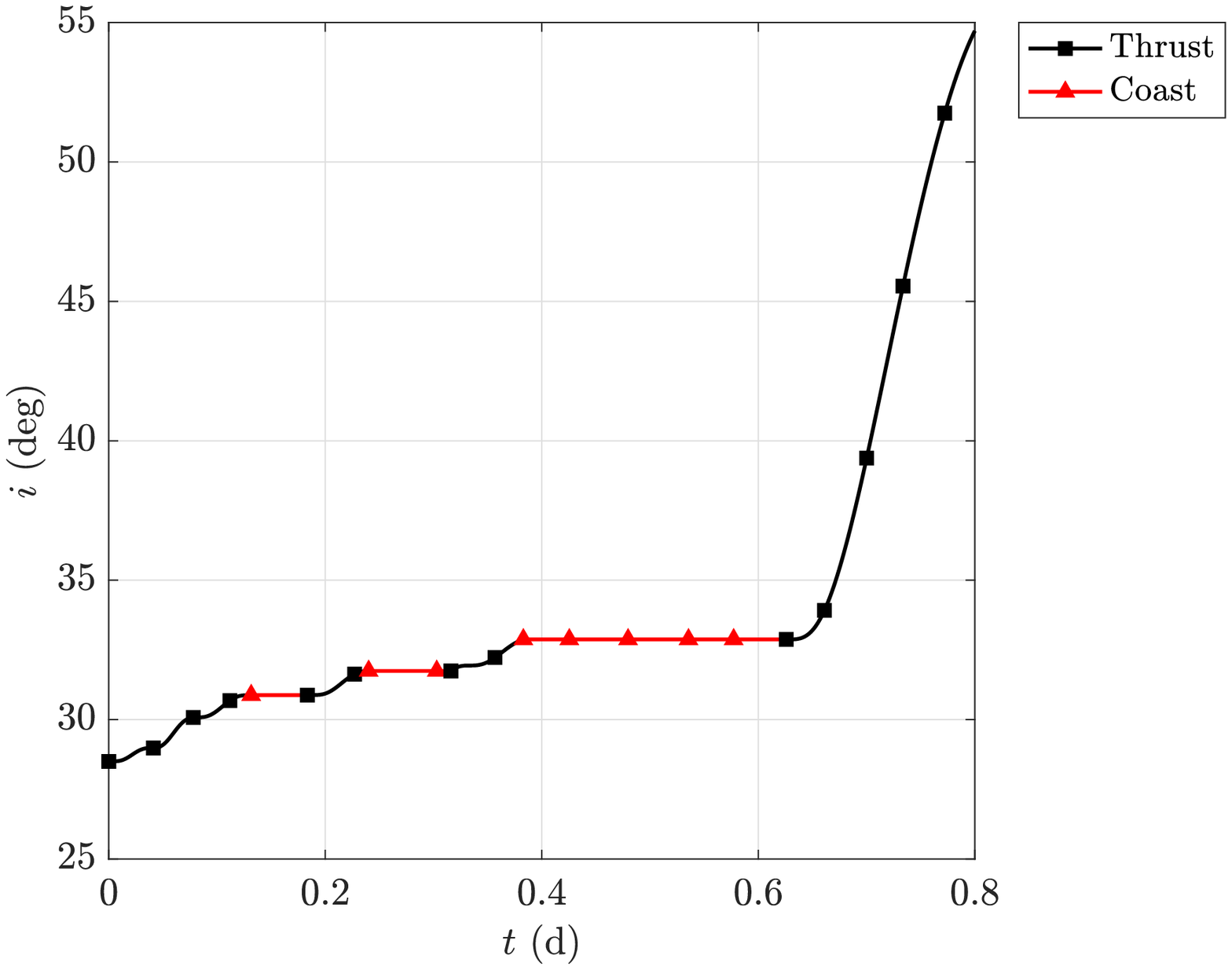}}
  \caption{Orbital elements along optimal trajectory for LEO-to-MEO transfer with $s_0 = 0.1~\textrm{m} \cdot \textrm{s}^{-2}$.\label{fig:MEO_Tm-1_orbital_elements}}
\end{figure}

Figure~\ref{fig:MEO_Tm-1_orbital_elements} shows the behavior of the orbital elements of the optimized trajectory of the spacecraft. The semi-major axis, $a$, increases throughout all of the thrust arcs from $7.003 \times 10^{6}~\textrm{m}$ to $2.6560 \times 10^{7}~\textrm{m}$, where the change becomes more rapid in the later thrust arcs. It is more fuel efficient to change the size of the orbit near periapsis because the velocity of the spacecraft is the fastest the spacecraft will travel on that specific orbit, meaning that the amount of fuel expended is less to achieve a higher velocity in the same direction than anywhere else on the orbit. Therefore, the thrust arcs happen near periapsis. The eccentricity, $e$, increases through the first three thrust arcs from $0$ to $0.3579$, then decreases to $0$ during the last thrust arc. The eccentricity rapidly changes during the first three thrust arcs because it is more fuel efficient to increase the size of the orbit first and get further from the Earth, then followed by changing the inclination of the orbit. The last thrust arc then creates a circular orbit. The inclination, $i$, increases slowly from $28.5~\textrm{deg}$ to $32.8764~\textrm{deg}$ during the first three thrust arcs, and then rapidly increases to $54.7~\textrm{deg}$ during the final thrust arc. The inclination changes by $4.3764~\textrm{deg}$ during the first three thrust arcs and by $21.8236~\textrm{deg}$ during the final thrust arc. The inclination changes much more significantly during the final thrust arc than in the first three thrust arcs because the spacecraft is farther away from the Earth, therefore the velocity of the spacecraft is smaller. Consequently, the maneuver is more fuel efficient because inclination changes require a change in the direction of velocity. Therefore, when the velocity is smaller the maneuver will require less fuel to be expended, so the inclination changes more during the final thrust arc. 

\begin{figure}[htb]
\centering
\includegraphics[width=3in]{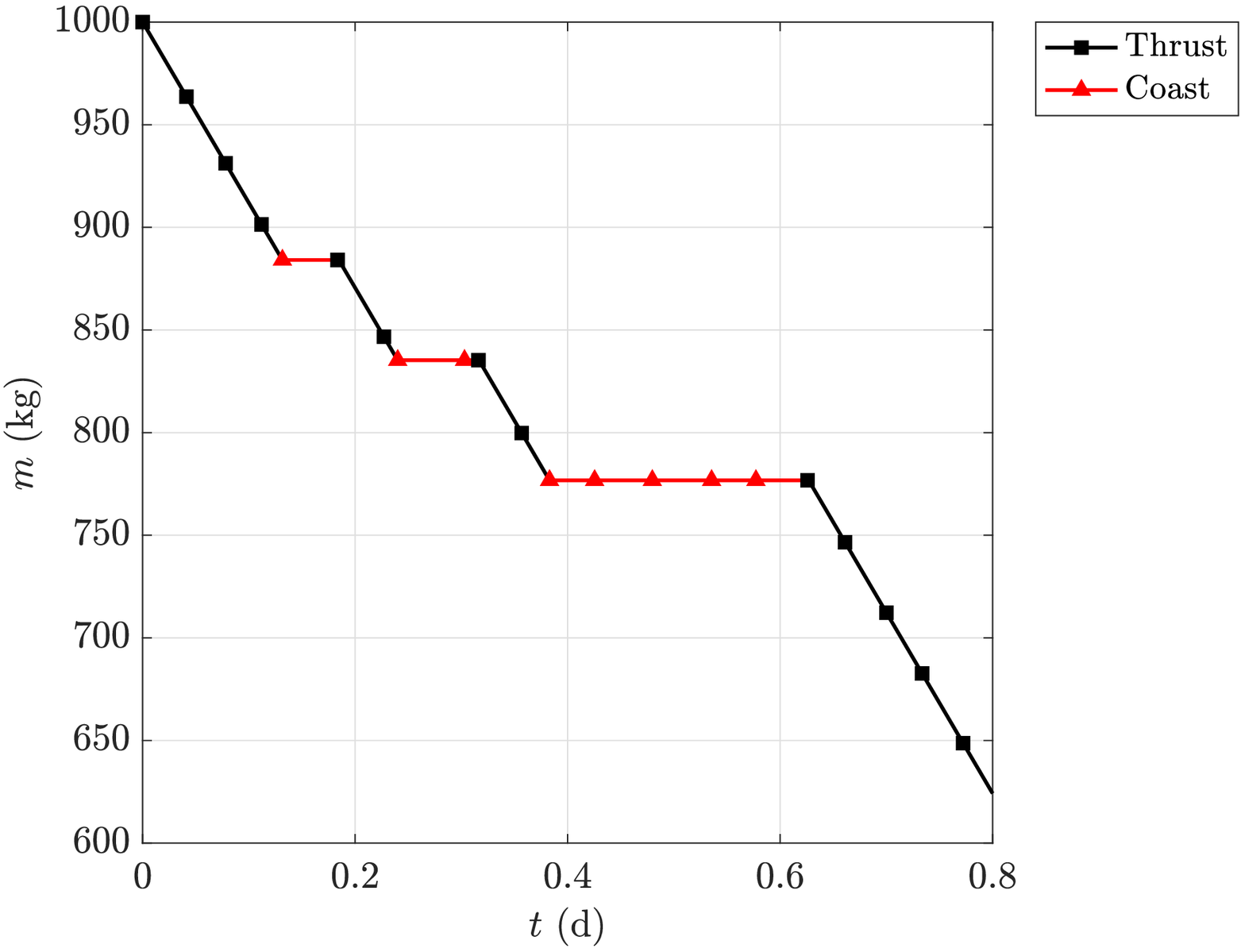}
\caption{Mass of spacecraft along optimal trajectory for LEO-to-MEO transfer with $s_0 = 0.1~\textrm{m} \cdot \textrm{s}^{-2}$.\label{fig:MEO_Tm-1_mvst}}
\end{figure}

Figure~\ref{fig:MEO_Tm-1_mvst} shows the mass of the spacecraft throughout the fuel-optimized trajectory. The mass decreases steadily throughout the four thrust arcs from $1000~\textrm{kg}$ to $624.2352~\textrm{kg}$, which is the optimized final mass. For Case 5 of the LEO-to-MEO transfer, the amount of fuel expended is $375.7648~\textrm{kg}$.

\begin{figure}[h]
  \centering
  \subfloat[Thrust magnitude.\label{fig:MEO_Tm-1_TTvst}]{\includegraphics[width=3in]{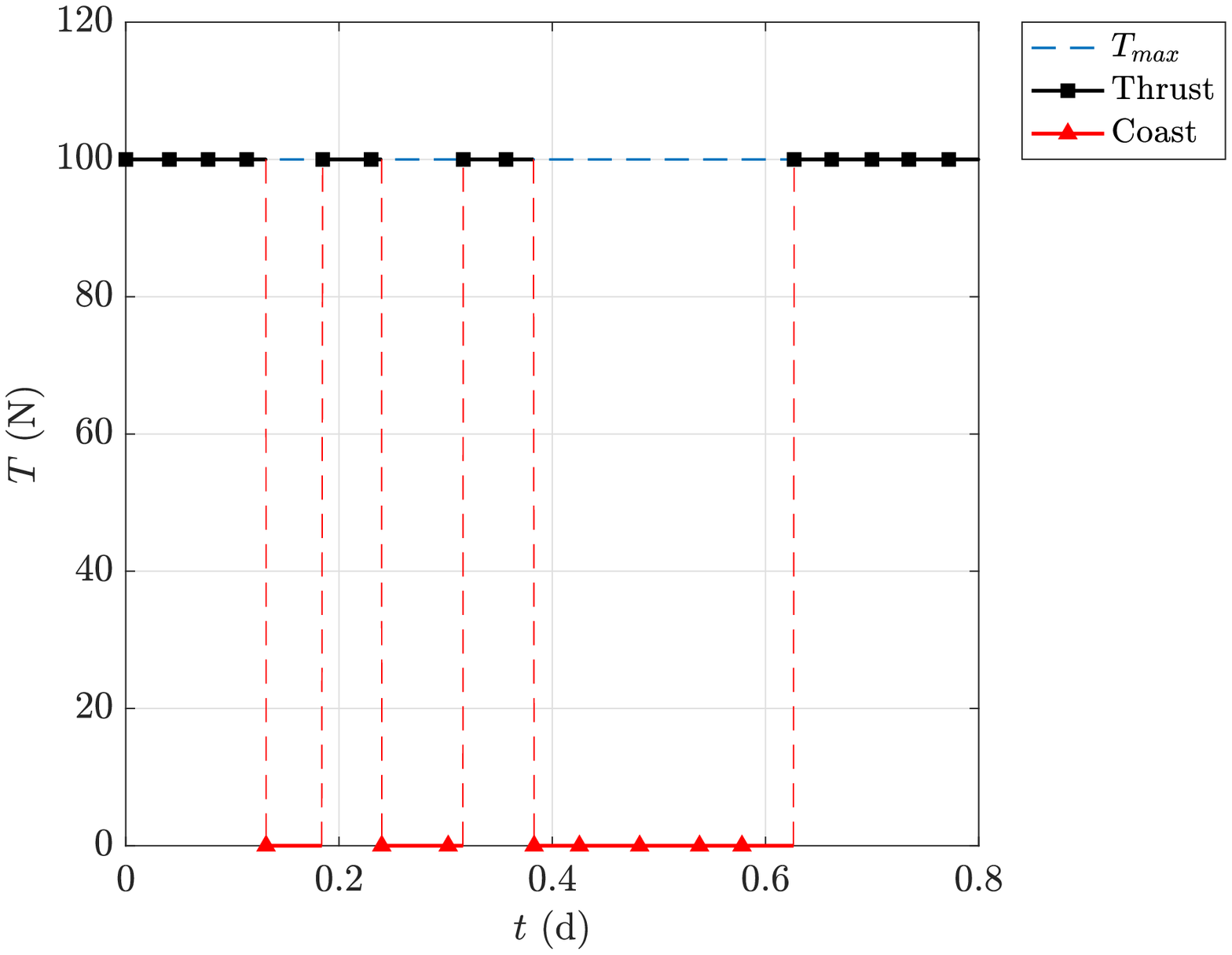}}~~\subfloat[Thrust direction components. \label{fig:MEO_Tm-1_uvst}]{\includegraphics[width=3in]{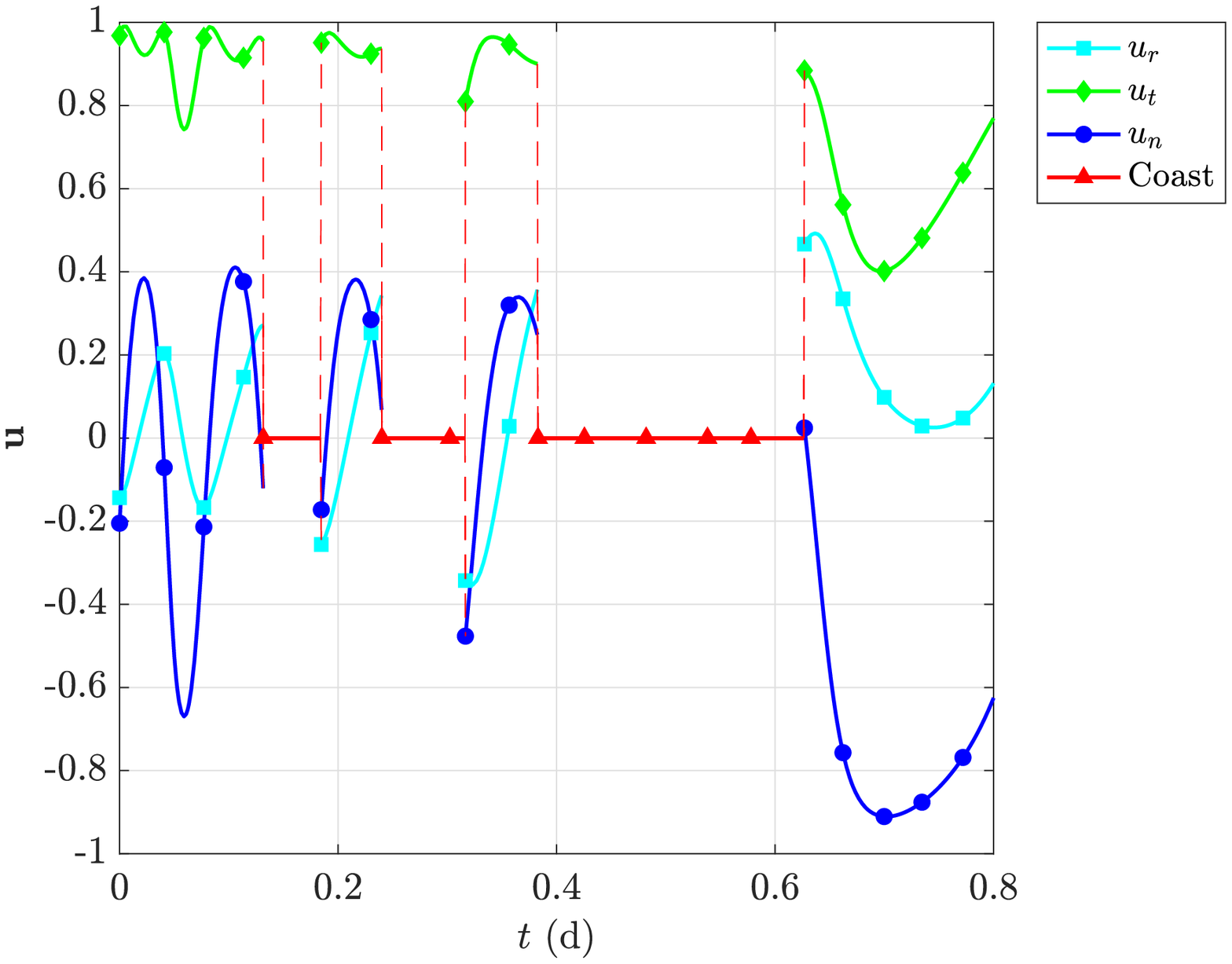}}
  \caption{Optimal control for LEO-to-MEO transfer with $s_0 = 0.1~\textrm{m} \cdot \textrm{s}^{-2}$.\label{fig:MEO_Tm-1_control}}
\end{figure}

Finally, Fig.~\ref{fig:MEO_Tm-1_control} shows the control components of the optimal trajectory, which are the thrust magnitude, $T$, and the thrust direction components, $\left(u_r,u_t,u_n\right)$. The thrust magnitude remains at the maximum thrust $T_{\max} = 100~\textrm{N}$ for the duration of the thrust arcs and $0~\textrm{N}$ for the duration of the coast arcs. There are four thrust arcs and three coast arcs in the optimized thrusting structure. Consequently, the thrust has six discontinuities and the structure of this solution is bang-bang. It is noted that the structure of the thrust is not assumed before solving the problem and the structure is detected using the BBSOC method~\cite{Pager2022}. The components of the thrust direction are only applicable when the thrust is non-zero (that is, the four thrust arcs), therefore the behavior will only be discussed for the thrust arcs because the components are set to zero during the coast arcs on the plot for clarity. The radial thrust direction component, $u_r$, increases from $-0.1435$ to $0.3569$ throughout the first three thrust arcs, where the direction oscillates about $0$ with increasing amplitudes, then decreases from $0.4665$ to $0.1319$ during the final thrust arc. The transverse thrust direction component, $u_t$, decreases from $0.9682$ to $0.9004$ throughout the first three thrust arcs, then decreases from $0.8842$ to $0.7693$ in the final thrust arc. The normal thrust direction component, $u_n$, increases from $-0.2048$ to $0.2488$ throughout the first three thrust arcs, where the direction oscillates, then decreases from $0.0242$ to $-0.6251$ during the final thrust arc. This behavior demonstrates that during the first three thrust arcs the majority of the thrust is in the transverse direction in order to increase the size of the orbit from $7.003 \times 10^{6}~\textrm{m}$ to $1.5862 \times 10^{7}~\textrm{m}$ and that during the final thrust arc the majority of the thrust is in the transverse and normal directions to increase the size of the orbit from $1.5862 \times 10^{7}~\textrm{m}$ to $2.6560 \times 10^{7}~\textrm{m}$ and increase the inclination of the orbit from $32.8764~\textrm{deg}$ to $54.7~\textrm{deg}$.

\subsubsection{Key Features of Optimized LEO-to-HEO Transfers}\label{sect: Study 2: LEO-to-HEO Key Features of Optimized Solutions}

\begin{figure}[htb]
\centering
\includegraphics[width=3in]{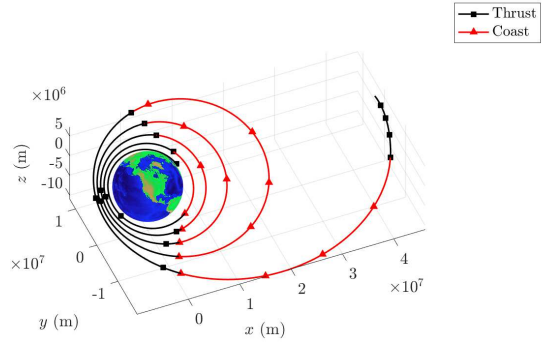}
\caption{Optimal three-dimensional spacecraft trajectory for LEO-to-HEO transfer with $s_0 = 0.1~\textrm{m} \cdot \textrm{s}^{-2}$.\label{fig:HEO_Tm-1_3D_Trajectory}}
\end{figure}

Figure~\ref{fig:HEO_Tm-1_3D_Trajectory} shows the optimized three-dimensional trajectory of the LEO-to-HEO transfer, where the modified equinoctial elements were converted into scaled Cartesian coordinates~\cite{MEE}. It is seen that the spacecraft begins in a low-Earth orbit and terminates in a high-Earth orbit that corresponds to the orbital elements in Table~\ref{tab:OrbitalElements}. The optimal trajectory of the spacecraft consists of $4.9570$ orbital revolutions around the Earth with a final mass of $657.2695~\textrm{kg}$, six thrust arcs, a total time thrusting of $9.2494~\textrm{h}$, and a total impulse of $4115.6~\textrm{m} \cdot \textrm{s}^{-1}$. It is noted that all of the thrust arcs, except for the final thrust arc, occur near periapsis of the transfer orbit. 

\begin{figure}[h]
  \centering
  \subfloat[Semi-major axis.\label{fig:HEO_Tm-1_avst}]{\includegraphics[width=3in]{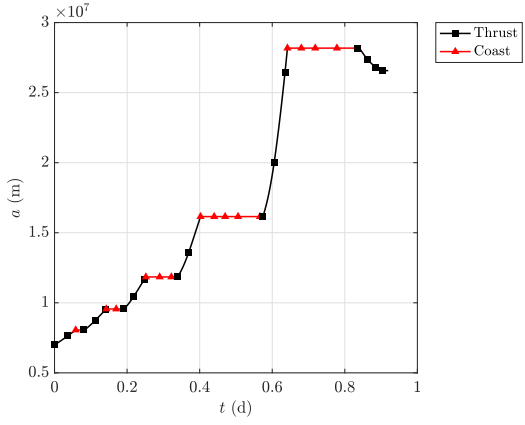}}~~\subfloat[Eccentricity. \label{fig:HEO_Tm-1_evst}]{\includegraphics[width=3in]{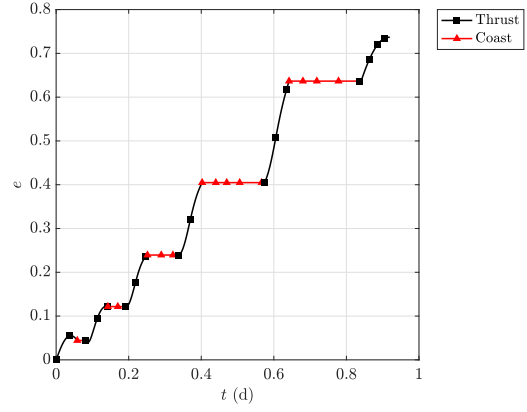}}

  \subfloat[Inclination. \label{fig:HEO_Tm-1_ivst}]{\includegraphics[width=3in]{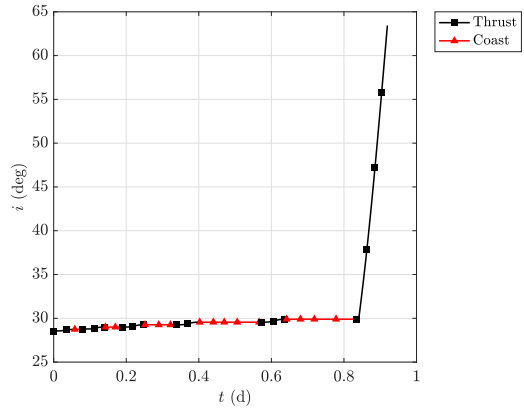}}
  \caption{Orbital elements along optimal trajectory for LEO-to-HEO transfer with $s_0 = 0.1~\textrm{m} \cdot \textrm{s}^{-2}$.\label{fig:HEO_Tm-1_orbital_elements}}
\end{figure}

Figure~\ref{fig:HEO_Tm-1_orbital_elements} shows the behavior of the orbital elements of the optimized trajectory of the spacecraft. The semi-major axis, $a$, increases throughout the first five thrust arcs from $7.003 \times 10^{6}~\textrm{m}$ to $2.8177 \times 10^{7}~\textrm{m}$, and then decreases to $2.6578 \times 10^{7}~\textrm{m}$ during the final thrust arc. It is more fuel efficient to change the size of the orbit near periapsis because the velocity of the spacecraft is the fastest the spacecraft will travel on that specific orbit, meaning that the amount of fuel expended is less to achieve a higher velocity in the same direction than anywhere else on the orbit. Therefore, the size of the orbit is changed the most during the fifth thrust arc. The eccentricity, $e$, increases throughout all of the thrust arcs from $0$ to $0.73646$. The eccentricity rapidly changes during the first five thrust arcs because it is more fuel efficient to change the size of the orbit first and get far from the Earth, then followed by changing the inclination of the orbit. The inclination, $i$, slowly increases from $28.5~\textrm{deg}$ to $29.8969~\textrm{deg}$ during the first five thrust arcs, and then rapidly increases to $63.435~\textrm{deg}$ during the final thrust arc. The inclination changes by $1.3969~\textrm{deg}$ during the first five thrust arcs and by $33.5381~\textrm{deg}$ during the final thrust arc. The inclination changes much more significantly during the final thrust arc than in the first five thrust arcs because the spacecraft is farther away from the Earth, therefore the velocity of the spacecraft is smaller. Consequently, the maneuver is more fuel efficient because inclination changes require a change in the direction of velocity. Therefore, when the velocity is smaller the maneuver will require less fuel to be expended, so the inclination changes more during the final thrust arc. 

\begin{figure}[htb]
\centering
\includegraphics[width=3in]{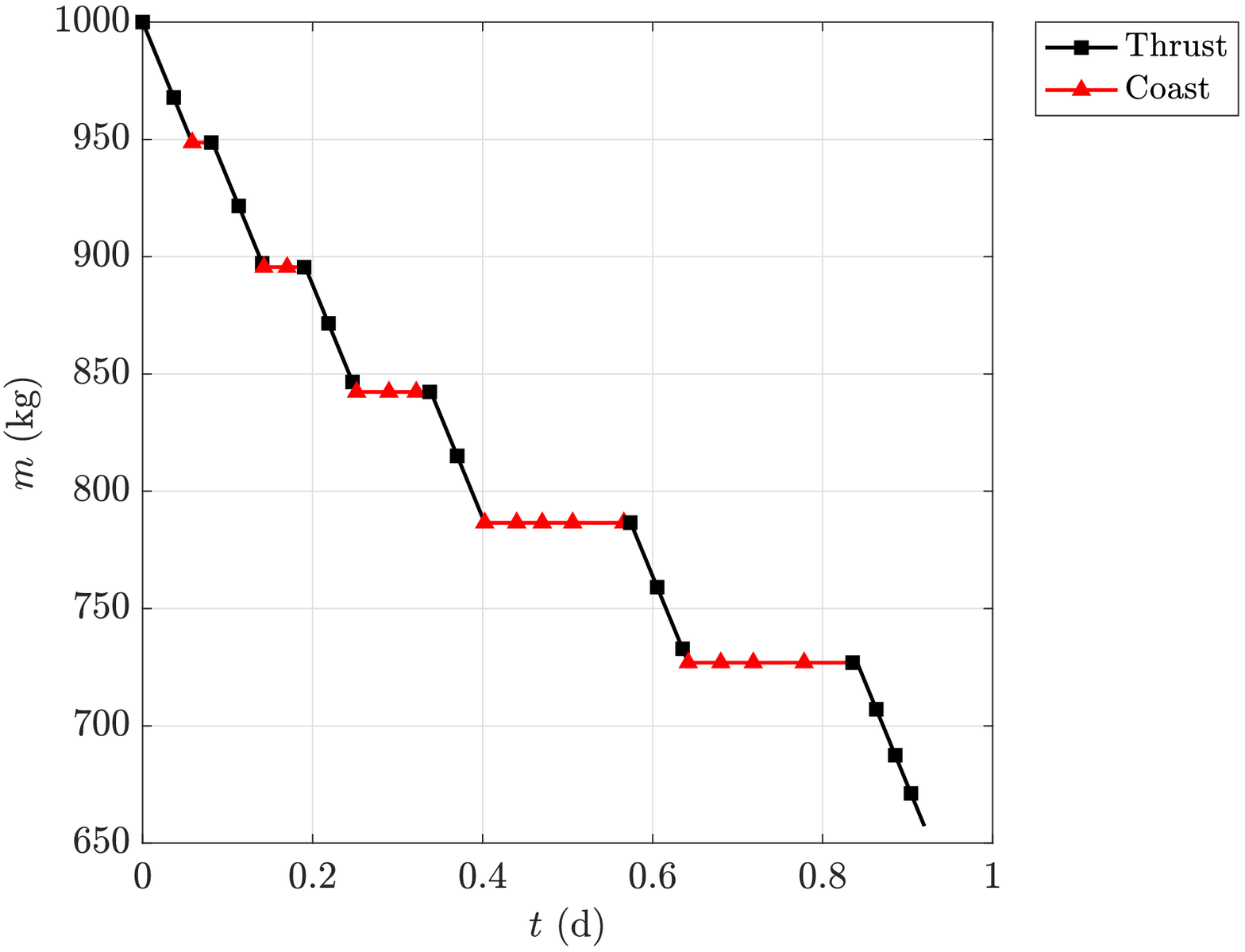}
\caption{Mass of spacecraft along optimal trajectory for LEO-to-HEO transfer with $s_0 = 0.1~\textrm{m} \cdot \textrm{s}^{-2}$.\label{fig:HEO_Tm-1_mvst}}
\end{figure}

Figure~\ref{fig:HEO_Tm-1_mvst} shows the mass of the spacecraft throughout the fuel-optimized trajectory. The mass decreases steadily throughout the six thrust arcs from $1000~\textrm{kg}$ to $657.2695~\textrm{kg}$, which is the optimized final mass. For Case 5 of the LEO-to-HEO transfer, the amount of fuel expended is $342.7305~\textrm{kg}$.

\begin{figure}[h]
  \centering
  \subfloat[Thrust magnitude.\label{fig:HEO_Tm-1_TTvst}]{\includegraphics[width=3in]{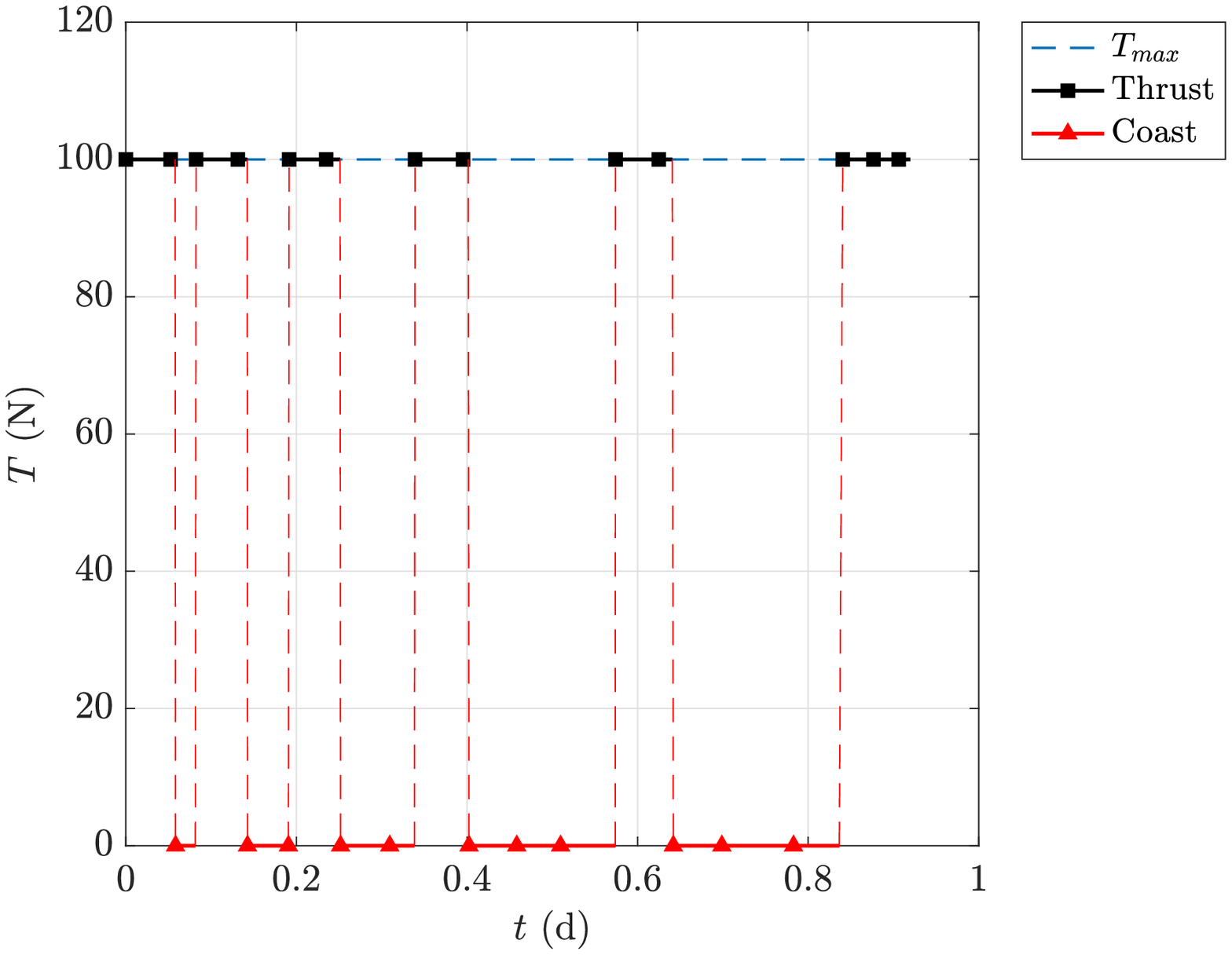}}~~\subfloat[Thrust direction components. \label{fig:HEO_Tm-1_uvst}]{\includegraphics[width=3in]{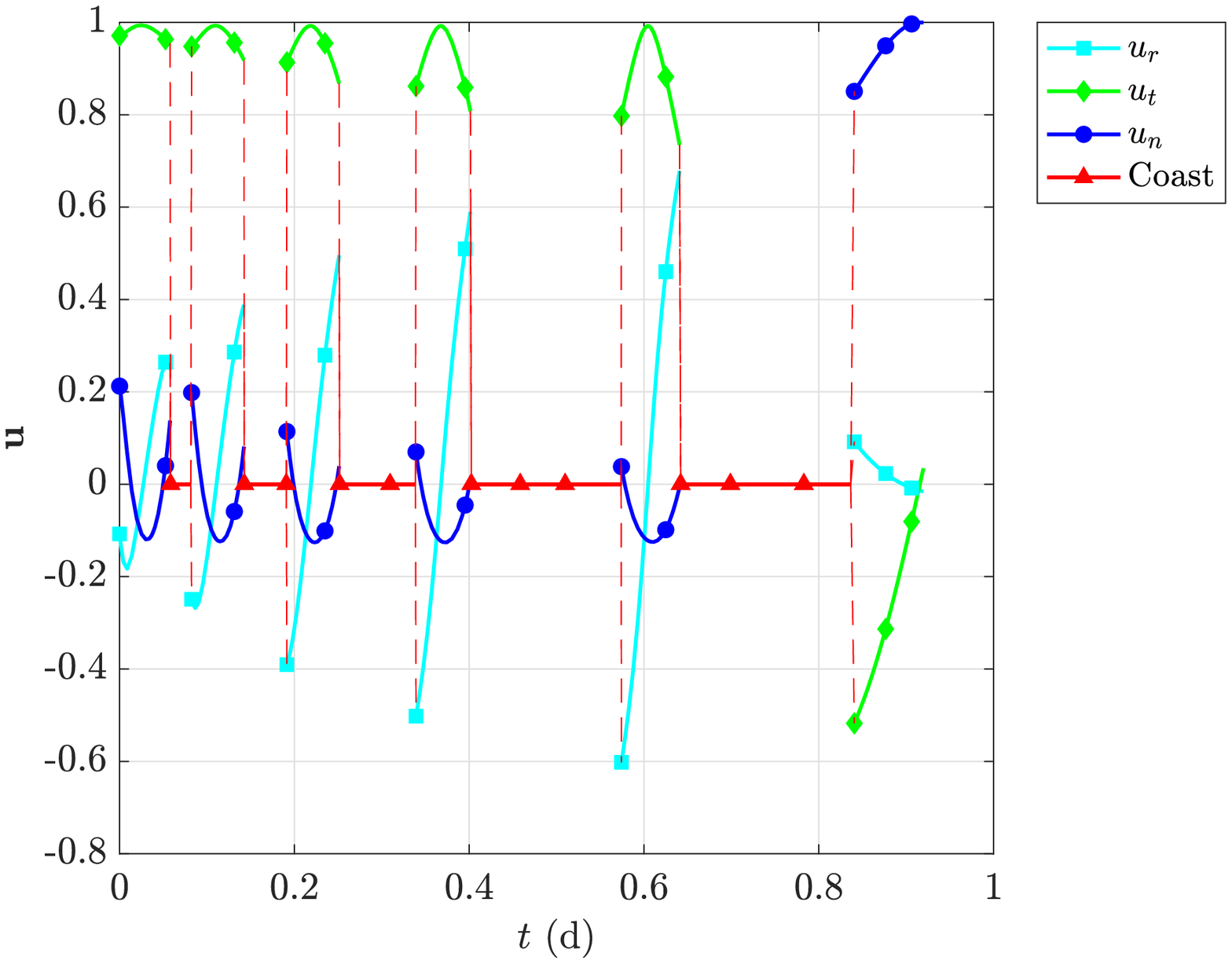}}
  \caption{Optimal control for LEO-to-HEO transfer with $s_0 = 0.1~\textrm{m} \cdot \textrm{s}^{-2}$.\label{fig:HEO_Tm-1_control}}
\end{figure}

Finally, Fig.~\ref{fig:HEO_Tm-1_control} shows the control components of the optimal trajectory, which are the thrust magnitude, $T$, and the thrust direction components, $\left(u_r,u_t,u_n\right)$. The thrust magnitude remains at the maximum thrust $T_{\max} = 100~\textrm{N}$ for the duration of the thrust arcs and $0~\textrm{N}$ for the duration of the coast arcs. There are six thrust arcs and five coast arcs in the optimized thrusting structure. Consequently, the thrust has ten discontinuities and the structure of this solution is bang-bang. It is noted that the structure of the thrust is not assumed before solving the problem and the structure is detected using the BBSOC method~\cite{Pager2022}. The components of the thrust direction are only applicable when the thrust is non-zero (that is, the six thrust arcs), therefore the behavior will only be discussed for the thrust arcs because the components are set to $0$ during the coast arcs on the plot for clarity. The radial thrust direction component, $u_r$, increases from $-0.1074$ to $0.6782$ throughout the first five thrust arcs, where the direction oscillates about $0$ with increasing amplitudes, then decreases from $0.0920$ to $-0.0157$ during the final thrust arc. The transverse thrust direction component, $u_t$, decreases from $0.9712$ to $0.7348$ throughout the first five thrust arcs, where the direction changes in an inverted parabolic pattern during the thrust arcs, then increases from $-0.5176$ to $0.0346$ in the final thrust arc. The normal thrust direction component, $u_n$, decreases from $0.2125$ to $-0.0062$ throughout the first five thrust arcs, where the direction changes in a parabolic pattern during the thrust arcs, then increases from $0.8506$ to $0.9997$ during the final thrust arc. This behavior demonstrates that during the first five thrust arcs the majority of the thrust is in the transverse direction in order to increase the size of the orbit from $7.003 \times 10^{6}~\textrm{m}$ to $2.8177 \times 10^{7}~\textrm{m}$ and that during the final thrust arc the majority of the thrust is in the normal direction to increase the inclination of the orbit from $29.8969~\textrm{deg}$ to $63.435~\textrm{deg}$. 

\subsubsection{Key Features of Optimized LEO-to-GEO Transfers}\label{sect: Study 3: LEO-to-HEO Key Features of Optimized Solutions}

\begin{figure}[htb]
\centering
\includegraphics[width=3in]{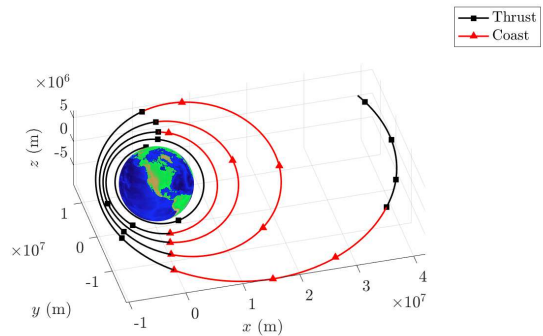}
\caption{Optimal three-dimensional spacecraft trajectory for LEO-to-GEO transfer with $s_0 = 0.1~\textrm{m} \cdot \textrm{s}^{-2}$.\label{fig:GEO_Tm-1_3D_Trajectory}}
\end{figure}

Figure~\ref{fig:GEO_Tm-1_3D_Trajectory} shows the optimized three-dimensional trajectory of the LEO-to-GEO transfer, where the modified equinoctial elements were converted into scaled Cartesian coordinates~\cite{MEE}. It is seen that the spacecraft begins in a low-Earth orbit and terminates in a geostationary orbit that corresponds to the orbital elements in Table~\ref{tab:OrbitalElements}. The optimal trajectory of the spacecraft consists of $4.8044$ orbital revolutions around the Earth with a final mass of $619.0090~\textrm{kg}$, five thrust arcs, a total time thrusting of $10.3158~\textrm{h}$, and a total impulse of $4703.6~\textrm{m} \cdot \textrm{s}^{-1}$. It is noted that all of the thrust arcs, except for the final, occur near the periapsis of the orbit. 

\begin{figure}[h]
  \centering
  \subfloat[Semi-major axis.\label{fig:GEO_Tm-1_avst}]{\includegraphics[width=3in]{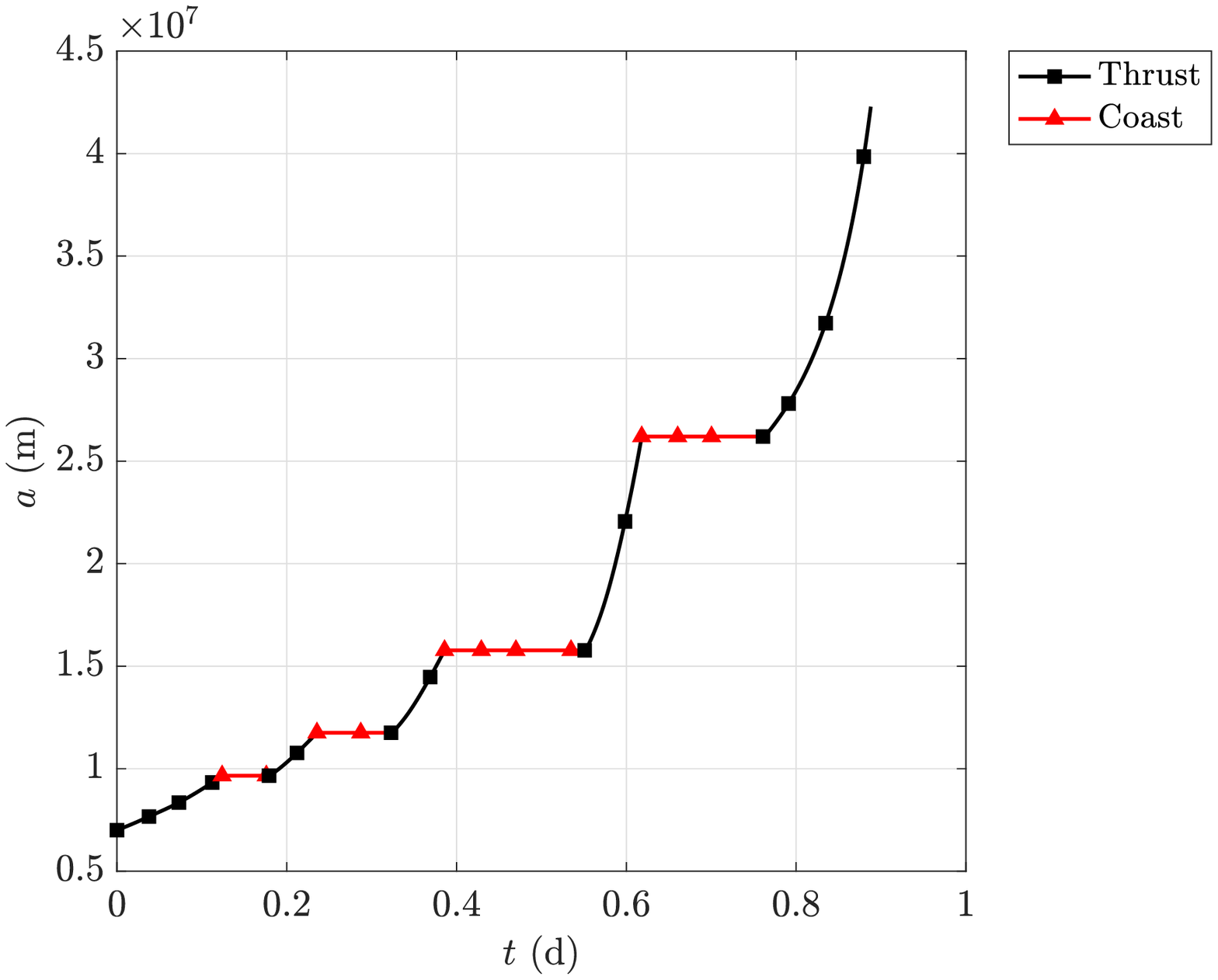}}~~\subfloat[Eccentricity. \label{fig:GEO_Tm-1_evst}]{\includegraphics[width=3in]{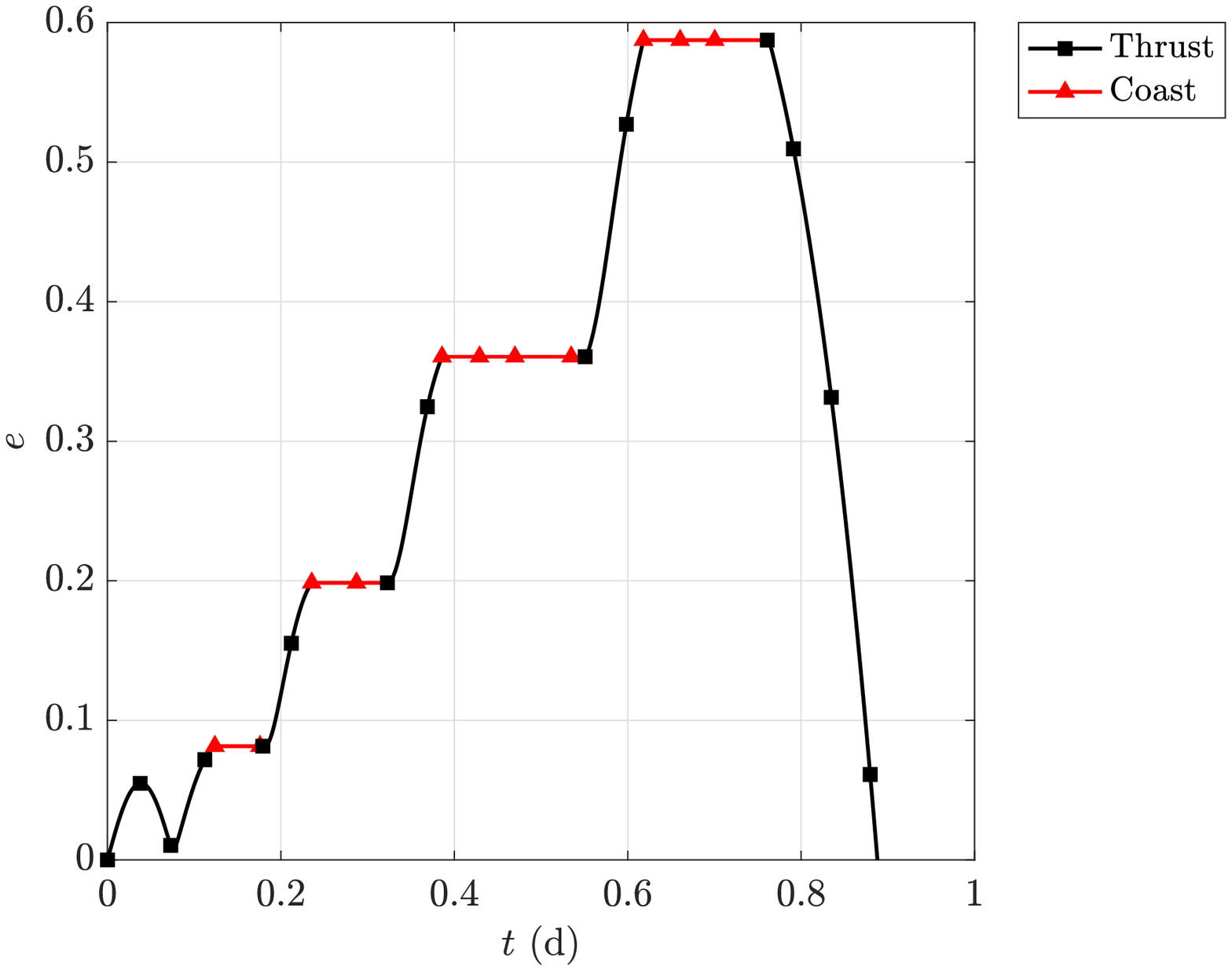}}

  \subfloat[Inclination. \label{fig:GEO_Tm-1_ivst}]{\includegraphics[width=3in]{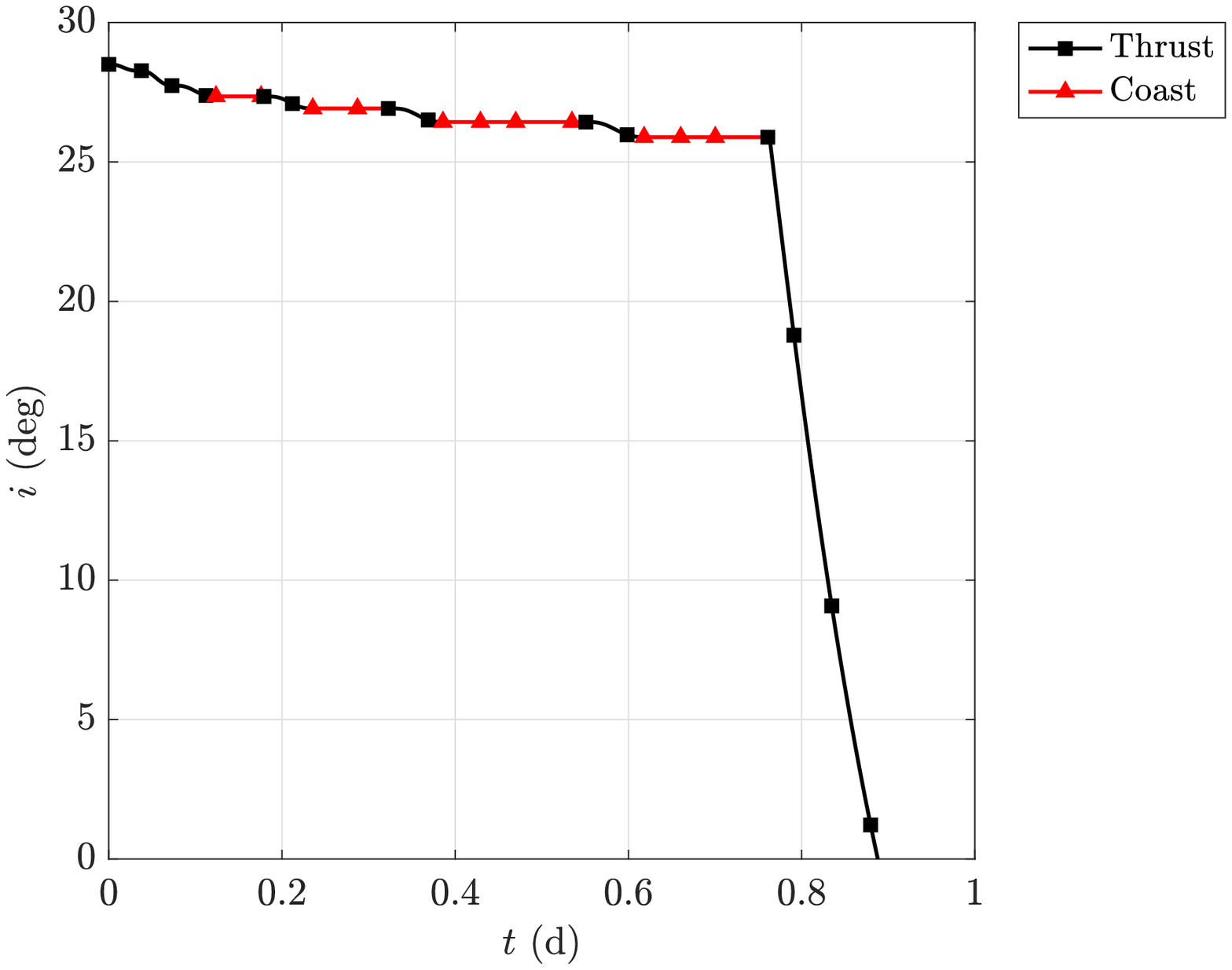}}
  \caption{Orbital elements along optimal trajectory for LEO-to-GEO transfer with $s_0 = 0.1~\textrm{m} \cdot \textrm{s}^{-2}$.\label{fig:GEO_Tm-1_orbital_elements}}
\end{figure}

Figure~\ref{fig:GEO_Tm-1_orbital_elements} shows the behavior of the orbital elements of the optimized trajectory of the spacecraft. The semi-major axis, $a$, increases throughout all of the thrust arcs from $7.003 \times 10^{6}~\textrm{m}$ to $4.2287 \times 10^{7}~\textrm{m}$, where the change becomes more rapid in the later thrust arcs. It is more fuel efficient to change the size of the orbit near periapsis because the velocity of the spacecraft is the fastest the spacecraft will travel on that specific orbit, meaning that the amount of fuel expended is less to achieve a higher velocity in the same direction than anywhere else on the orbit. Therefore, the thrust arcs happen near periapsis, except for the final thrust arc. The eccentricity, $e$, increases through the first four thrust arcs from $0$ to $0.5874$, then decreases to $0$ during the last thrust arc. The eccentricity rapidly changes during the first four thrust arcs because it is more fuel efficient to increase the size of the orbit first and get far from the Earth, then followed by changing the inclination of the orbit. The last thrust arc then creates a circular orbit. The inclination, $i$, slowly decreases from $28.5~\textrm{deg}$ to $25.8893~\textrm{deg}$ during the first four thrust arcs, and then rapidly decreases to $0~\textrm{deg}$ during the final thrust arc. The inclination changes by $2.6107~\textrm{deg}$ during the first four thrust arcs and by $25.8893~\textrm{deg}$ during the final thrust arc. The inclination changes much more significantly during the final thrust arc than in the first four thrust arcs because the spacecraft is farther away from the Earth, therefore the velocity of the spacecraft is smaller. Consequently, the maneuver is more fuel efficient because inclination changes require a change in the direction of velocity. Therefore, when the velocity is smaller the maneuver will require less fuel to be expended, so the inclination changes more during the final thrust arc. 

\begin{figure}[htb]
\centering
\includegraphics[width=3in]{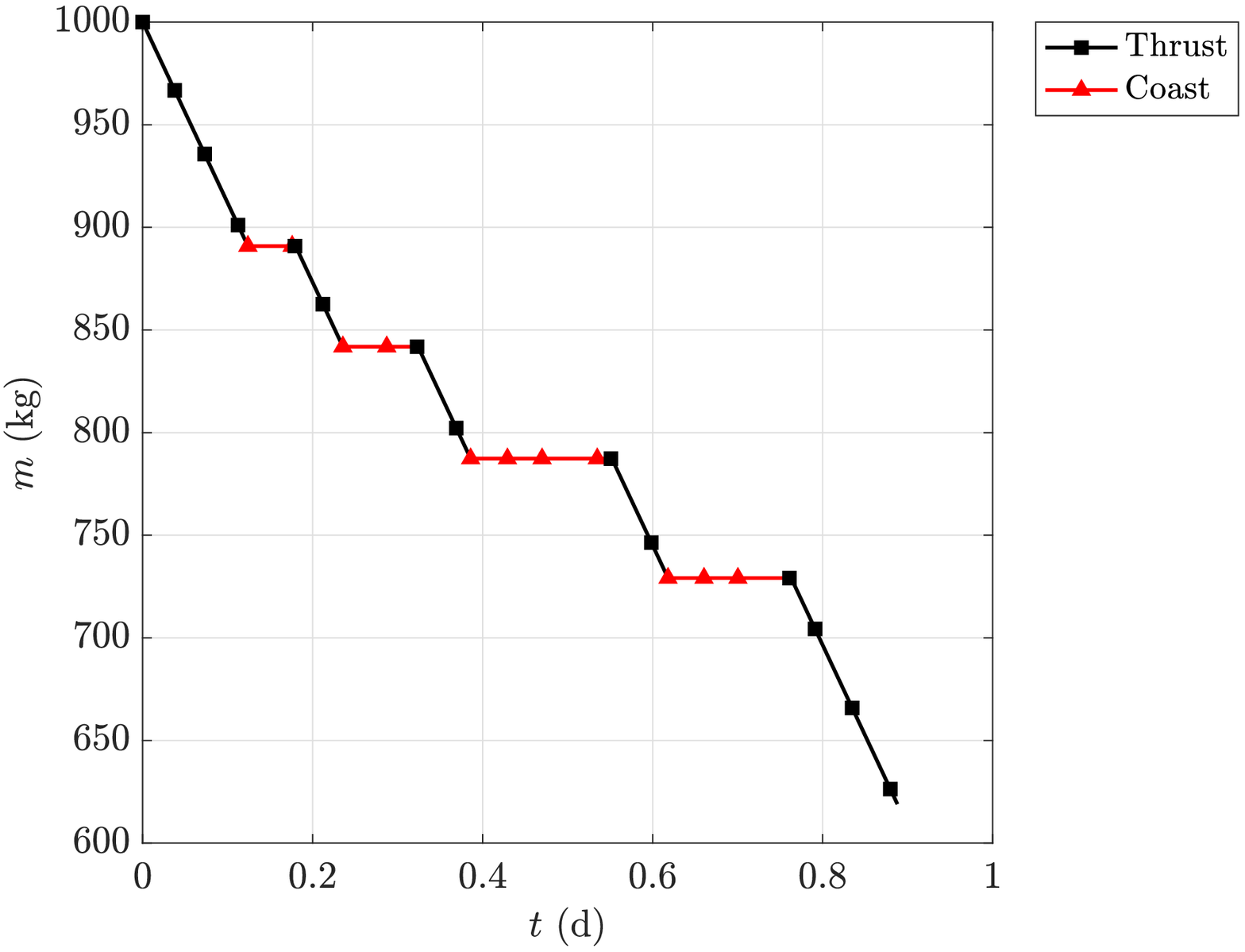}
\caption{Mass of spacecraft along optimal trajectory for LEO-to-GEO transfer with $s_0 = 0.1~\textrm{m} \cdot \textrm{s}^{-2}$.\label{fig:GEO_Tm-1_mvst}}
\end{figure}

Figure~\ref{fig:GEO_Tm-1_mvst} shows the mass of the spacecraft throughout the fuel-optimized trajectory. The mass decreases steadily throughout the five thrust arcs from $1000~\textrm{kg}$ to $619.0090~\textrm{kg}$, which is the optimized final mass. For Case 5 of the LEO-to-GEO transfer, the amount of fuel expended is $380.9910~\textrm{kg}$.

\begin{figure}[h]
  \centering
  \subfloat[Thrust magnitude.\label{fig:GEO_Tm-1_TTvst}]{\includegraphics[width=3in]{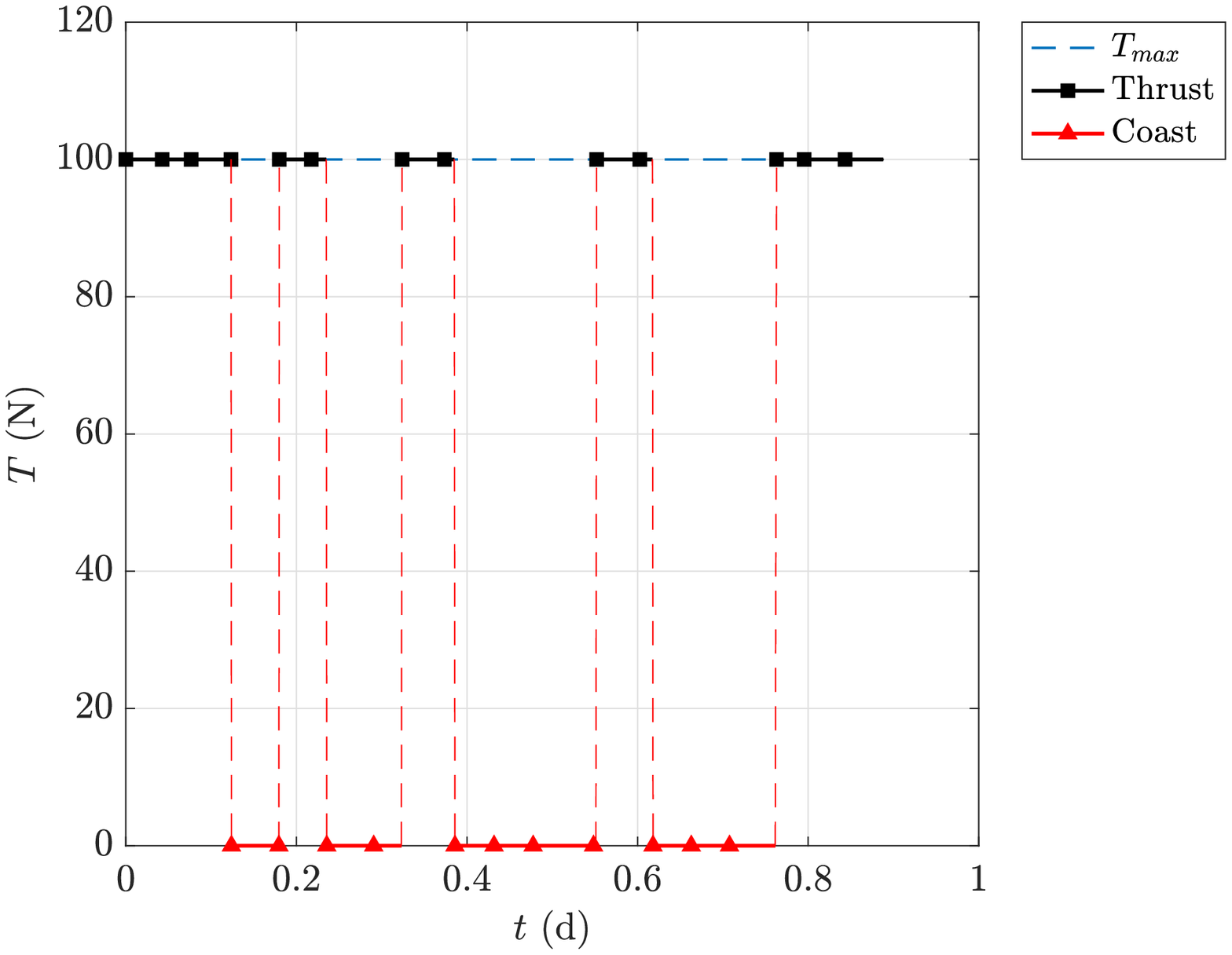}}~~\subfloat[Thrust direction components. \label{fig:GEO_Tm-1_uvst}]{\includegraphics[width=3in]{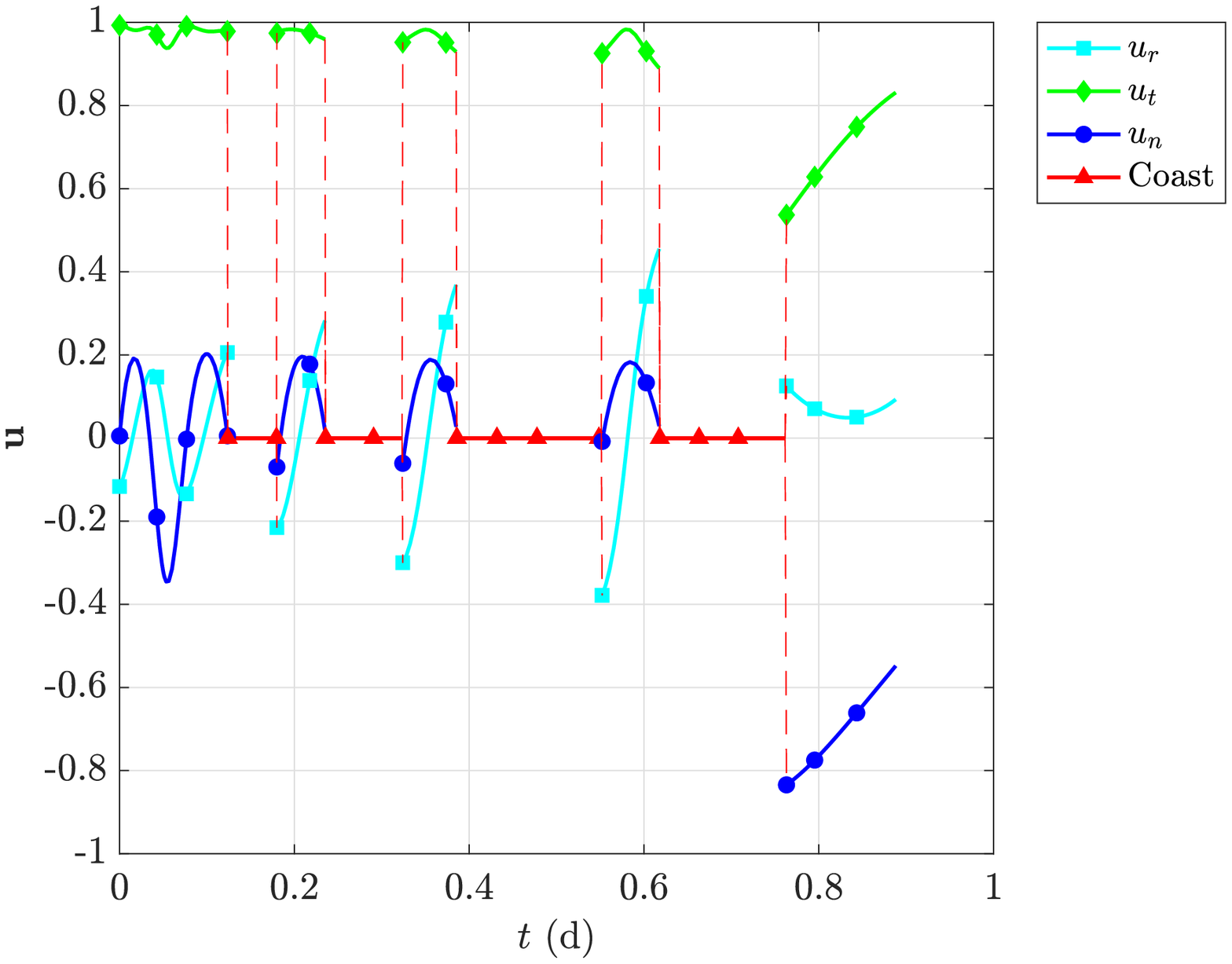}}
  \caption{Optimal control for LEO-to-GEO transfer with $s_0 = 0.1~\textrm{m} \cdot \textrm{s}^{-2}$.\label{fig:GEO_Tm-1_control}}
\end{figure}

Finally, Fig.~\ref{fig:GEO_Tm-1_control} shows the control components of the optimal trajectory, which are the thrust magnitude, $T$, and the thrust direction components, $\left(u_r,u_t,u_n\right)$. The thrust magnitude remains at the maximum thrust $T_{\max} = 100~\textrm{N}$ for the duration of the thrust arcs and $0~\textrm{N}$ for the duration of the coast arcs. There are five thrust arcs and four coast arcs in the optimized thrusting structure. Consequently, the thrust has eight discontinuities and the structure of this solution is bang-bang. It is noted that the structure of the thrust is not assumed before solving the problem and the structure is detected using the BBSOC method~\cite{Pager2022}. The components of the thrust direction are only applicable when the thrust is non-zero (that is, the five thrust arcs), therefore the behavior will only be discussed for the thrust arcs because the components are set to $0$ during the coast arcs on the plot for clarity. The radial thrust direction component, $u_r$, increases from $-0.1167$ to $0.4551$ throughout the first four thrust arcs, where the direction oscillates about $0$ with increasing amplitudes, then decreases from $0.1253$ to $0.0929$ during the final thrust arc. The transverse thrust direction component, $u_t$, decreases from $0.9932$ to $0.8900$ throughout the first four thrust arcs, then increases from $0.5368$ to $0.8313$ in the final thrust arc. The normal thrust direction component, $u_n$, increases from $0.0051$ to $0.0288$ throughout the first three thrust arcs, where the direction oscillates, then increases from $-0.8343$ to $-0.5481$ during the final thrust arc. This behavior demonstrates that during the first four thrust arcs the majority of the thrust is in the transverse direction in order to increase the size of the orbit from $7.003 \times 10^{6}~\textrm{m}$ to $2.6198 \times 10^{7}~\textrm{m}$ and that during the final thrust arc the majority of the thrust is in the transverse and normal directions to increase the size of the orbit from $2.6198 \times 10^{7}~\textrm{m}$ to $4.2287 \times 10^{7}~\textrm{m}$ and decrease the inclination of the orbit from $25.8893~\textrm{deg}$ to $0~\textrm{deg}$.

\section{Conclusions}\label{sect:Conclusions}

This research performed a numerical optimization study of various minimum-fuel Earth-based orbital trajectory transfers. A method to solve bang-bang optimal control problems using multi-domain Legendre-Gauss-Radau quadrature collocation was implemented. Three orbital transfers were studied, Study 1: LEO-to-MEO, Study 2: LEO-to-HEO, and Study 3: LEO-to-GEO, where each study was solved for seven cases of maximum allowable thrust acceleration. The time, $t$, state components, $\left(p,f,g,h,k,L,m\right)$, and control components, $\left(T,u_r,u_t,u_n\right)$ were optimized in this one-phase problem. The gravitational parameter of Earth was scaled to be unity in order to improve numerical performance. Two initial guess generation schemes were utilized depending on the type of solution that was produced, either a partial or multiple orbital revolution solution. A key feature of this research is that the thrusting structure was not assumed a priori, therefore the optimizer determined the number of switch points. The minimum-fuel results for each study all follow the same pattern: as the maximum allowable thrust acceleration decreases, the final mass decreases, the total time thrusting increases, the number of orbital revolutions increases, the number of thrust arcs increases, and the total impulse increases. In addition, in all studies Cases 1-4 were categorized as partial orbital revolution solutions and Cases 5-7 were categorized as multiple orbital revolution solutions. Furthermore, all partial orbital revolution solutions had two thrust arcs. For the cases of each study, the final mass decreases because the total time thrusting increases with decreasing maximum allowable thrust acceleration, which leads to more fuel being expended to achieve the same terminal conditions. The key features of the optimized solutions were shown only for Case 5, $s_0 = 0.1~\textrm{m} \cdot \textrm{s}^{-2}$. For Study 1: LEO-to-MEO, the optimal trajectory of the spacecraft consists of $4.9579$ orbital revolutions around the Earth with a final mass of $624.2352~\textrm{kg}$, four thrust arcs, a total time thrusting of $10.2132~\textrm{h}$, and a total impulse of $4621.2~\textrm{m} \cdot \textrm{s}^{-1}$. For Study 2: LEO-to-HEO, the optimal trajectory of the spacecraft consists of $4.9570$ orbital revolutions around the Earth with a final mass of $657.2695~\textrm{kg}$, six thrust arcs, a total time thrusting of $9.2494~\textrm{h}$, and a total impulse of $4115.6~\textrm{m} \cdot \textrm{s}^{-1}$. For Study 3: LEO-to-GEO, the optimal trajectory of the spacecraft consists of $4.8044$ orbital revolutions around the Earth with a final mass of $619.0090~\textrm{kg}$, five thrust arcs, a total time thrusting of $10.3158~\textrm{h}$, and a total impulse of $4703.6~\textrm{m} \cdot \textrm{s}^{-1}$. When comparing total impulse values of this research and that performed in Ref.~\cite{Herman2002} for the shared maximum allowable thrust acceleration cases, as the thrusting structure more closely resembles the assumed thrusting structure of burn-coast-burn in Ref.~\cite{Herman2002}, the total impulse obtained in this study matches more closely the total impulse obtained in Ref.~\cite{Herman2002}. On the other hand, it is seen that, when these thrusting structures are dissimilar, the total impulse obtained in this study is significantly lower than that obtained in Ref.~\cite{Herman2002}.

\section*{Acknowledgments}\label{sect:Acknowledgments}

The authors gratefully acknowledge support for this research from the U.S.~National Science Foundation under grant CMMI-2031213 and the National Aeronautics and Space Administration under grant NNX15AI10H through the University of Central Florida NASA Space Grant Consortium and Space Florida.

\renewcommand{\baselinestretch}{1}
\normalsize
\normalfont
\bibliographystyle{aiaa}

\end{document}